 \documentclass[11pt,twoside]{article}
 \usepackage[dvips]{graphics}
 \usepackage{amsmath,graphicx,amsfonts,amsrefs,amssymb}
 \usepackage{setspace, enumitem, mathtools, float, yfonts}
 \usepackage{graphicx}
 \usepackage[table]{xcolor}
 
 \DeclarePairedDelimiter{\abs}{\lvert}{\rvert}

 \def\qed{\hfill\rule{1ex}{1ex}\\}

\begin{document}
 
 \title{Spectral methods for capillary surfaces described by bounded generating curves}
 \author{ Ray Treinen \footnote{Department of Mathematics, Texas State University, 601 University Dr., San Marcos, TX 78666, rt30@txstate.edu}
 }
 \maketitle

\begin{abstract}
We consider  capillary surfaces that are constructed by bounded generating curves.  This class of  surfaces includes radially symmetric and lower dimensional fluid-fluid interfaces.  We use the arc-length representation of the differential equations for these surfaces to allow for vertical points and inflection points along the generating curve.  These considerations admit capillary tubes, sessile drops, and fluids in annular tubes as well as other examples. 

We present a pseudo-spectral method for approximating  solutions to the associated boundary value problems based on interpolation by Chebyshev polynomials.  This method is observably more stable than the traditional shooting method and it is computationally lean and fast.  The algorithm is also adaptive, but does not use the adaptive automation in Chebfun.
\end{abstract}


\section{Introduction}
\label{intro}

The equilibrium shape of liquids can be described by the interfaces between immiscible fluids, and those interfaces are commonly called capillary surfaces.  These  fluid shapes have been of interest since the ancient times \cite{Mccuan2009}, though the advent of the calculus, and soon thereafter that of the calculus of variations gave precision to these studies, and in 1806 Laplace gave, among other results, the first description of radially symmetric capillary surfaces that this author is aware of.  The study of capillary action with radial symmetry was the motivation for the numerical method we now know as Adams-Bashforth \cite{AdamsBashforth1883},  a classic multi-step ODE solver.  The dawn of the space age renewed interest in these capillary problems, as surface and wetting effects dominate fluid behavior in micro-gravity.  The monograph by Finn \cite{ecs}  is an excellent resource for those unfamiliar with capillarity.  Around the turn of the century interest swung towards problems involving multiple fluids in equilibrium, or configurations of fluids with floating solids \cite{BhatnagarFinn2006}, \cite{ENS2004}.
We will not  give a complete catalog of these results here.  We will note that, as is extremely common in mathematical physics, the numerical computation of solutions is far ahead of the theoretical treatment.  It is our goal here to provide a robust numerical solver for those capillary surfaces that are described by bounded generating curves which is both computationally lean and extremely precise.   This will provide a useful tool in the further study of configurations where at least one component  is a capillary surface.

Chebyshev polynomials were first described in 1854 \cite{Chebyshev1854}.  Their use in approximation theory has been extensive, and a catalog of the uses of these orthogonal polynomials will also not be given here.  We will, however, echo Trefethen's observation \cite{Trefethen2013} that Orszag's early work in Chebyshev spectral methods dates to the beginning of the 1970's in \cite{Orszag1971a} and \cite{Orszag1971b}. We will also be using these polynomials to interpolate solutions, though for us they will be solutions to the differential equations describing capillary surfaces, and this will be in the form of Chebyshev differentiation matrices used as building blocks for the differential operators involved.  For those unfamiliar with this technique or background, we recommend Trefethen \cite{Trefethen2000} and \cite{Trefethen2013},  though there are other fine books (e.g.  \cite{Boyd2001} and \cite{Fornberg1996}),.  Trefethen, Birkisson, and Driscoll \cite{TrefethenBirkissonDriscoll2018} is also quite relevant to our approach.  We will be using Matlab with Chebfun \cite{Chebfun} to facilitate some aspect of our work, though we will not be using all of the automation implemented in Chebfun.  We prefer to present enough details of our work so that an interested reader can implement these algorithms on platforms other than Matlab using Chebfun without an undue burden to develop the code.  With this aim for ease of cross-platform implementation, we use Chebfun to generate our differentiation matrices, and for plotting our generating curves with barycentric interpolation in a way that puts visible points at the Chebyshev points along the curves.

\begin{figure}[!t]
	\centering
	\scalebox{0.35}{\includegraphics{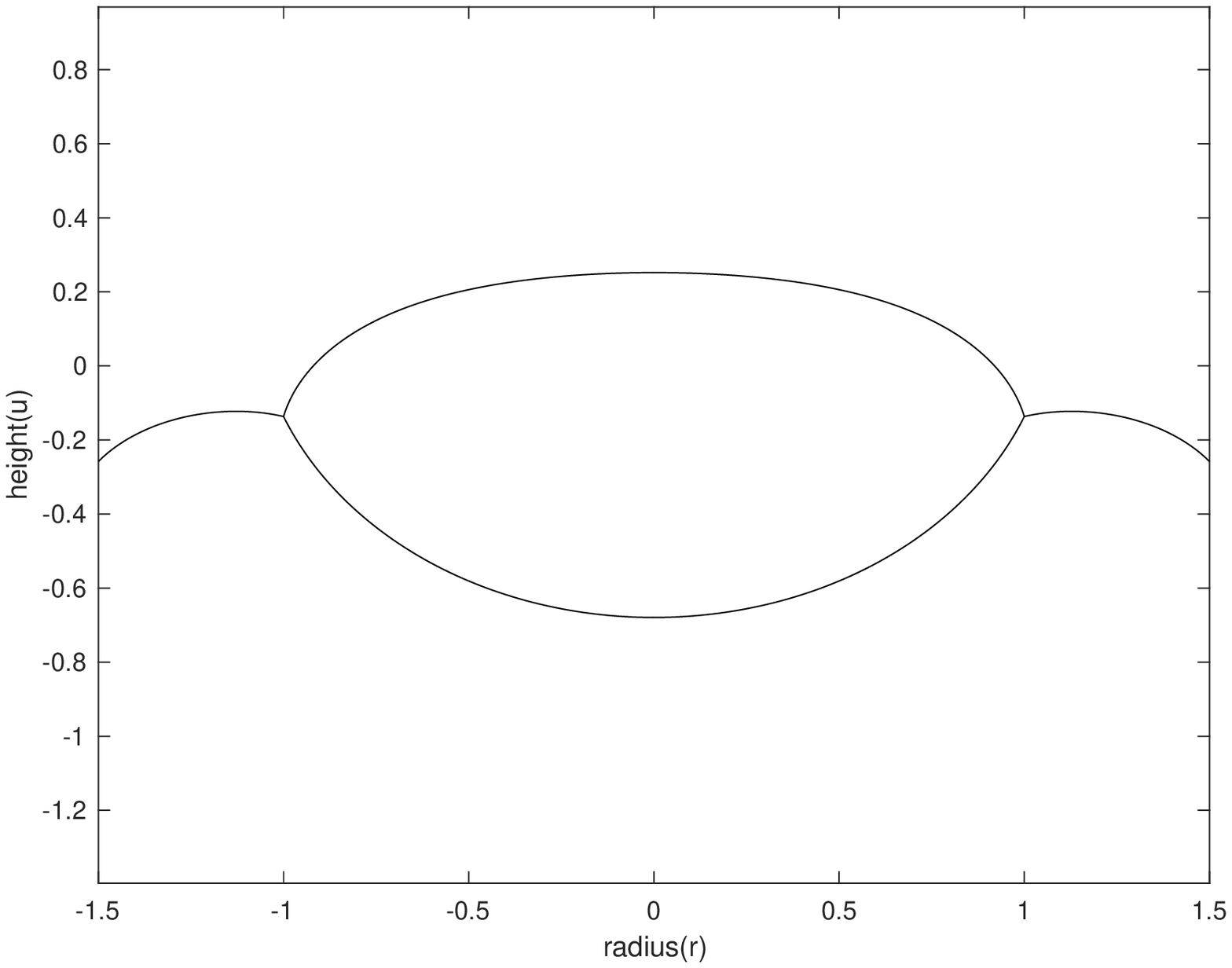}}
	\scalebox{0.35}{\includegraphics{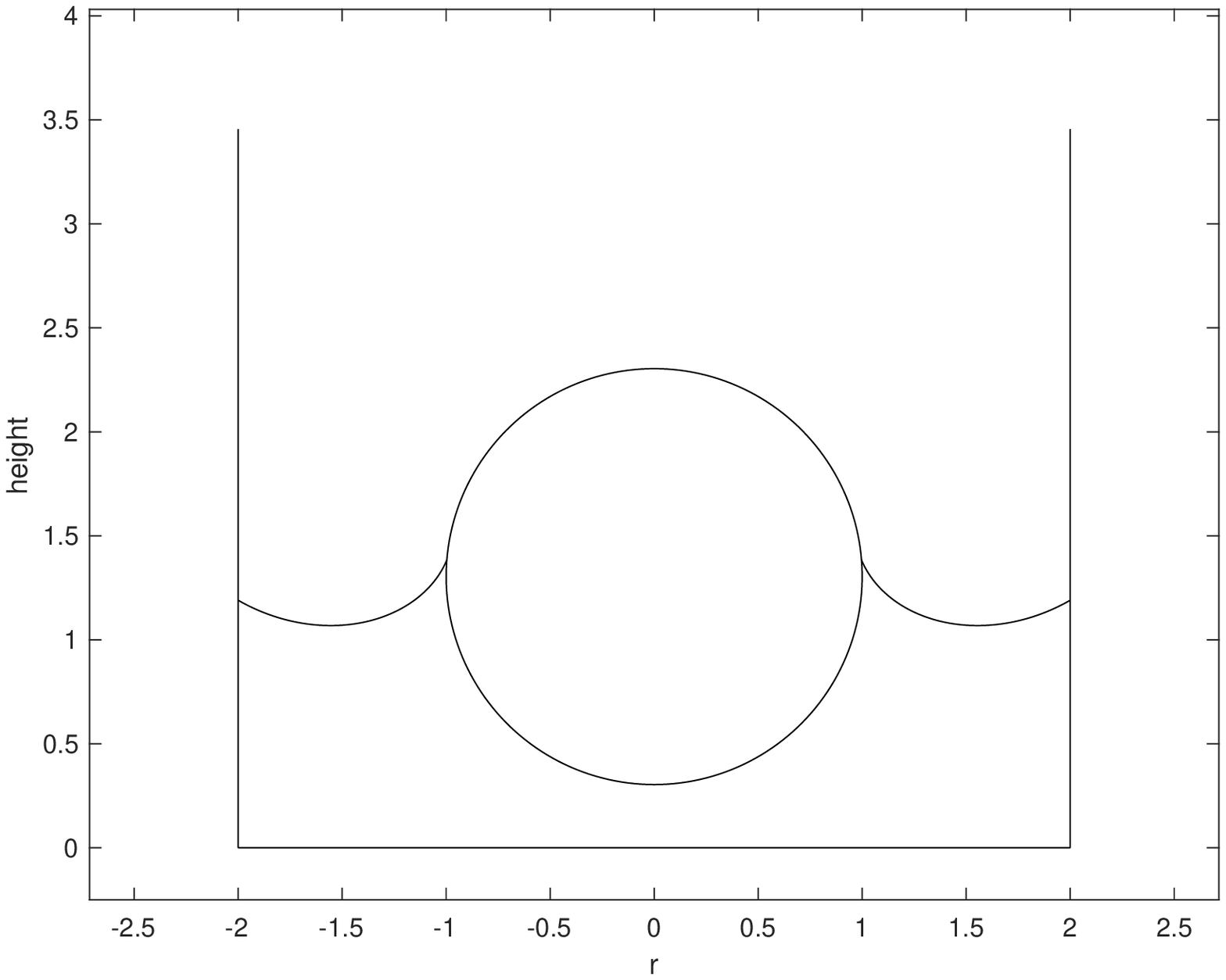}}
	\caption{A drop of liquid rests at the interface between water and air (L) and a ball floats at the surface of a liquid. (R).  Both examples are sections of radially symmetric configurations, and both use legacy code to generate the capillary surfaces shown.}
	\label{fig:examples2}
\end{figure}

We will focus on three prototype problems:
\begin{itemize}
	\item {\bf P1}
	Simply connected interfaces that are the image of a disk,
	\item {\bf P2}
	Doubly connected interfaces that are the image of an annulus, and
	\item {\bf P3}
	The lower-dimensional problem.
\end{itemize}
\begin{figure}[!b]
	\centering
	\scalebox{0.35}{\includegraphics{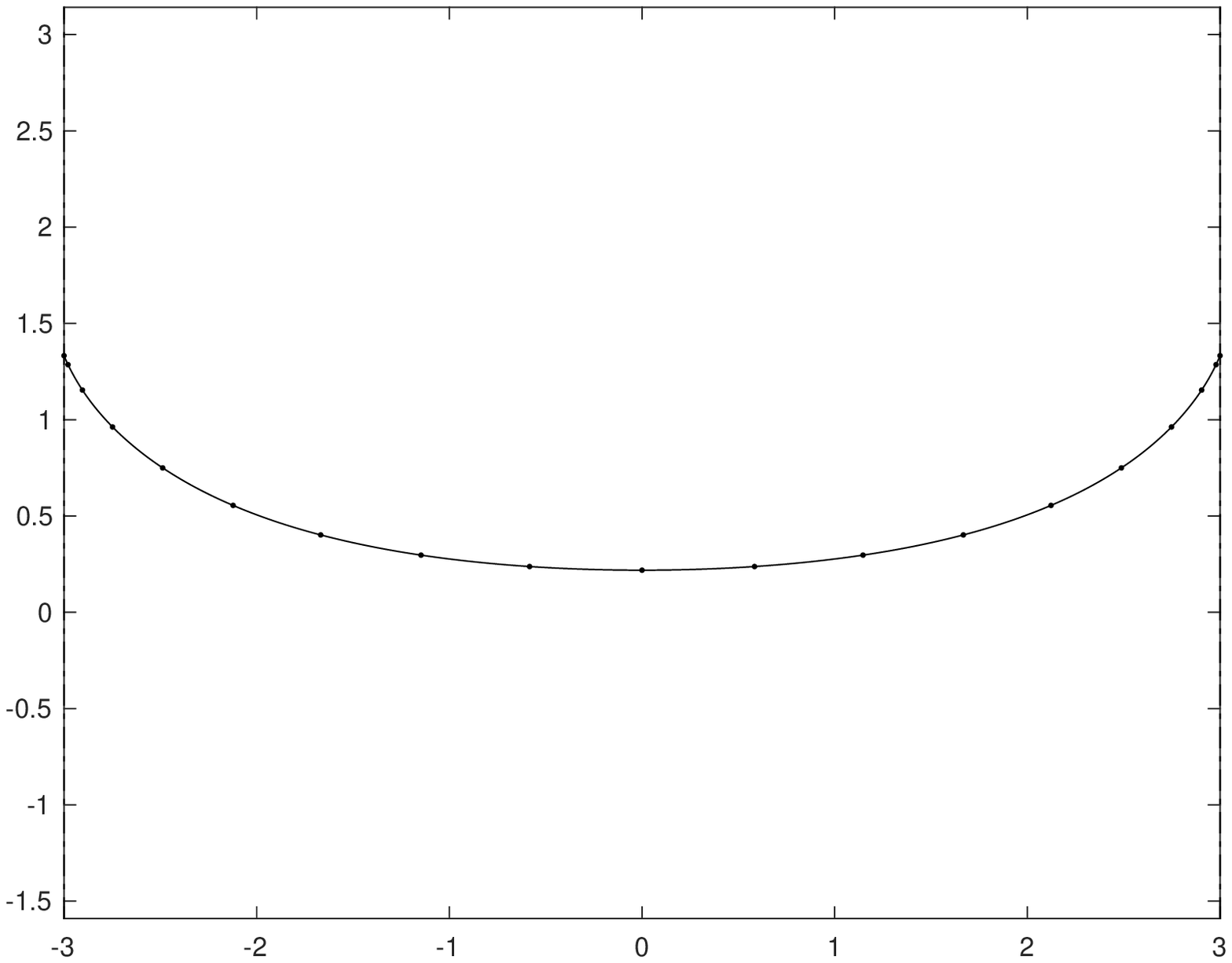}}
	\scalebox{0.35}{\includegraphics{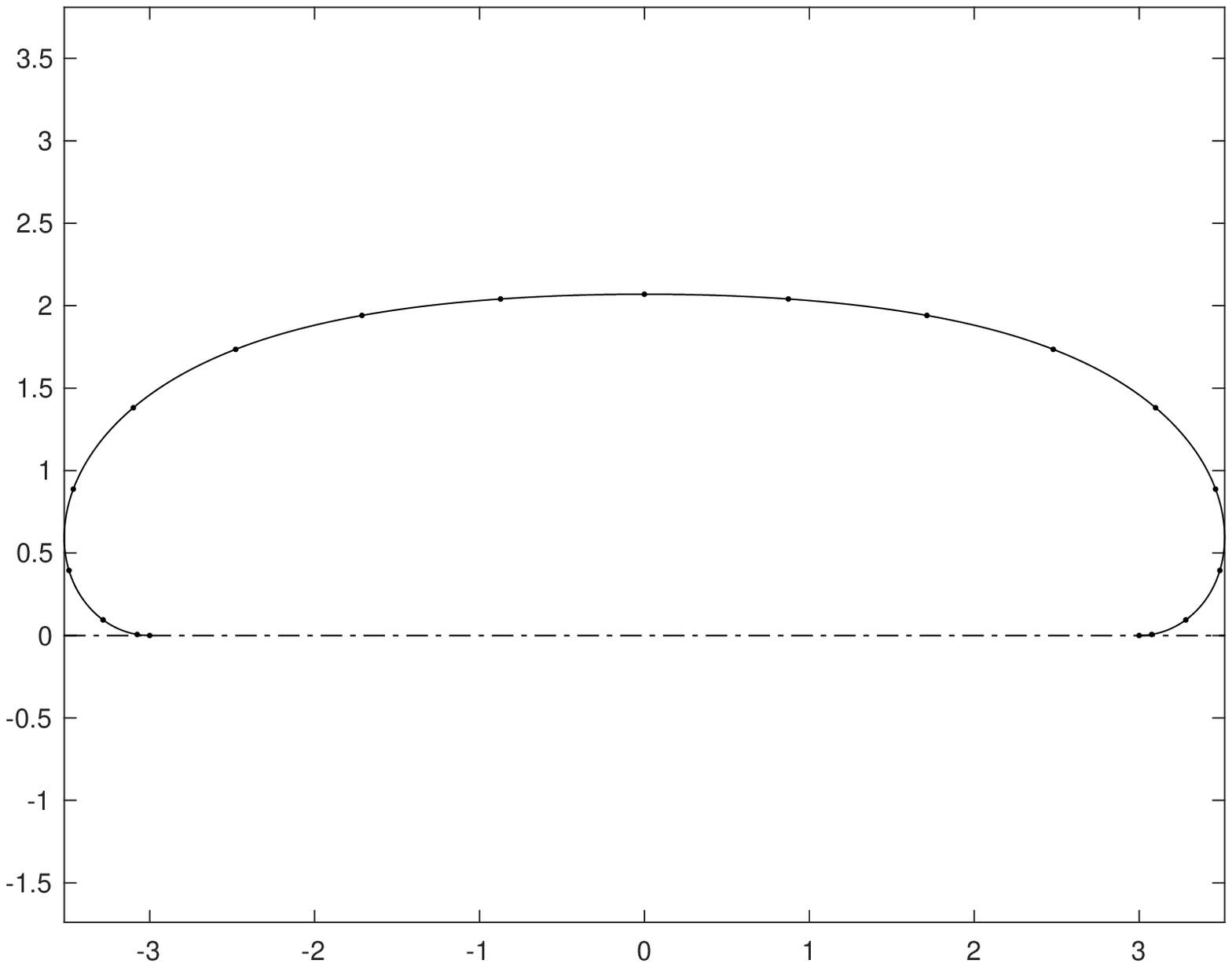}}
	\caption{A capillary tube with radius 3 and contact angle $\gamma = 3\pi/8$ (L) and a sessile drop with contact angle $\gamma = \pi$ (R).}
	\label{fig:examples}
\end{figure}
The main examples for {\bf P1} are the capillary tube and the sessile drop: see Figure~\ref{fig:examples}.  
The main examples for {\bf P2} are component interfaces for more complicated configurations such as a container filled with three immiscible fluids in equilibrium seen as a floating drop and a tube partially filled with a liquid where a ball is floating in the center of the tube, both shown in Figure~\ref{fig:examples2}.  The lower dimensional versions of these multiple-fluid configurations or liquid and floating solid configurations can provide insight when the full-dimensional problem remains largely unexplored.  For this reason, we include the details of {\bf P3}.  We will also use this setting to comment on  the new method developed here compared with the shooting method that has been used for {\bf P1}-{\bf P3} over the previous two decades.

The rest of this paper is organized as follows: in Section~\ref{Details} we give details for {\bf P1}-{\bf P3},   in Section~\ref{Cheby} we give cover the basic algorithm, in Sections~\ref{Study1}-\ref{Study3} we give studies of each treatment of {\bf P1}-{\bf P3} and provide details specific to each case.  Section ~\ref{conclusions} contains conclusions and a discussion of the limitations of the code.

We have released Matlab implementations of these algorithms under an open source license.  These files are available in a GitHub repository at \\
\texttt{https://github.com/raytreinen/Spetral-Method-Capillarity.git} \\
Finally, all of the testing for this paper was done on a 2018 MacBook Pro with a 2.6 GHz 6-Core Intel Core i7 processor and 32 GB 2400 MHz DDR4 RAM.

\section{Details for the three prototype problems}
\label{Details}

The mean curvature of a capillary surface is an affine function of its height  $u$ over some reference level.  We will restrict our attention to embedded solutions of the mean curvature equation
$$
nH = \kappa u - \lambda,
$$
where $nH$ is $n$ times the mean curvature of the surface, $\kappa = \rho g / \sigma > 0$ is the capillary constant with $\rho$ defined to be the difference in the densities of the fluids, $g$ the gravitational constant, and $\sigma$ the surface tension.  Here $\lambda$ is a Lagrange multiplier that is used for volume-constrained problems, and for the physical problems that we are considering, $n$ is 1 or 2, depending on the dimension of the ``surface''.  See \cite{ecs} and \cite{Wente2006}.  This equation is sometimes known as the Young-Laplace equation, or the capillary equation.

We first note that vertically translating a surface does not change its mean curvature, and so we may treat our problems without the Lagrange multiplier $\lambda$ and then vertically translate the surface to meet any volume constraints that may exist.  With the asumption that $\lambda = 0$, we next note that if $u$ is a solution, then $-u$ is also a solution.  Our next observation is that our mean-curvature formulation is not dependent on the surface being described as a function over a base domain, and this feature is often  found in treatments of these problems.  We wish to unify treatments of sessile drops with treatments of fluids in capillary tubes, and other potential component interfaces, and so we choose formulations that can include vertical points and inflection points on the capillary surface.  In case the reader is not familiar, we note that this strategy appears in \cite{ecs}, and we refer to that work for further reading.

In the present work we will be restricting our attention to surfaces that are described by a bounded generating curve.  These are most commonly  surfaces that are radially symmetric.  In some cases there are theorems that directly establish this symmetry, and in other cases the existence of symmetric solutions is found, and then the general comparison theorem for these surfaces can be applied to show uniqueness, and thus the symmetric solution is the unique solution.  There are configurations where symmetry is not yet established, and other configurations where the symmetric solutions are shown to have a larger potential energy than asymmetric solutions. (See \cite{ConcusFinn1991} and \cite{ConcusFinnWeislogel1999}.)  In any case, we find ongoing interest in the symmetric solutions of capillary surfaces.

With this radial symmetry, the capillary equation can be written as a system of three nonlinear ordinary differential  equations, parametrized by the arclength $s$:
\begin{eqnarray}
\frac{dr}{ds} &=& \cos\psi, \label{drds}\\
\frac{du}{ds} &=& \sin\psi, \label{duds}\\
\frac{d\psi}{ds} &=& \kappa u - \frac{\sin\psi}{r}, \label{dpsids}
\end{eqnarray}
where $r$ is the radius and $u$ is the height of the interface, and the inclination angle $\psi$ is measured from the corresponding generating curve which is described by $(r(s),u(s))$.

\subsection{ The problem {\bf P1} }
\label{subsecP1}

For symmetric and simply connected surfaces that are the image of a disk, we specify boundary conditions by the requirement that at some arclength $\ell>0$ the radius $r(\ell)$ meets a prescribed value $b>0$, and the inclination angle $\psi(\ell)$ meets a prescribed value $\psi_b \in [-\pi, \pi]$.  However, this value of the arclength $\ell$ is unknown, and this predisposes one to a shooting method to determine this value with other parameters to solve the problem.  We take an alternative approach by rescaling the problem.  We define $\tau = s/\ell$, or $s = \ell\tau$.  Then we define
\begin{eqnarray}
R(\tau) &:=& r(\ell\tau) = r(s),  \nonumber\\
U(\tau) &:=& u(\ell\tau) = u(s), \nonumber\\
\Psi(\tau) &:=& \psi(\ell\tau) = \psi(s). \nonumber 
\end{eqnarray}
This scaling is suitable for use with Chebyshev polynomials, as the natural domains of $R,U,$ and $\Psi$ are now $[-1,1]$.  Then, using the chain rule and multiplying each equation by $\ell$, \eqref{drds}-\eqref{dpsids} become
\begin{eqnarray}
R^\prime(\tau) - \ell\cos\Psi(\tau) &=& 0, \label{eqn:R}\\
U^\prime(\tau) - \ell\sin\Psi(\tau) &=& 0, \label{eqn:U}\\
\Psi^\prime(\tau) + \frac{\ell\sin\Psi(\tau)}{R(\tau)} - \kappa\ell U(\tau) &=& 0. \label{eqn:Psi}
\end{eqnarray}
If we define the column vector $\mathbf{v} = [R\,\, U\,\, \Psi\,  \ell]^T$, we can use \eqref{eqn:R}-\eqref{eqn:Psi} to define the nonlinear operator in the vector equation
\begin{equation}
	\tilde N(\mathbf{v}) = \mathbf{0}. \label{eqn:tN=0}
\end{equation}
We then use the boundary conditions
\begin{eqnarray}
R(1) - b &=& 0, \label{eqn:bcRb} \\
R(-1) + b &=& 0, \label{eqn:bcRmb} \\
\Psi(1) - \psi_b &=& 0, \label{eqn:bcPsib} \\
\Psi(-1) + \psi_b &=& 0, \label{eqn:bcmPsib}
\end{eqnarray}
so that we have some form of a two-point boundary value problem.  We append \eqref{eqn:tN=0} with these  boundary conditions to form the system 
\begin{equation}
{N}(\mathbf{v}) = \mathbf{0}.
\end{equation}
The resulting solution $(R,U) = (r,u)$ is not a generating curve of the surface, but it is a section.  The generating curve is given by $(R(\tau),U(\tau))$ for $0\leq\tau\leq 1$, and the generating curve has total arclength $\ell$. It is worth noting, that in our formulation there is no explicit restriction to curves that have a reflection symmetry about the vertical axis.  Of course, only those symmetric configurations are solutions, but the numerical solver will be free to choose asymmetric curves. One final remark on the general nature of this problem is that \eqref{dpsids} and \eqref{eqn:Psi} have singularities at $r = R = 0$.  It is known that this singularity is removable, and the solution is analytic \cite{ecs}.  We will discuss our numerical approach to this singularity when we return to the results for this problem.  For the time being, we present an alternative form of \eqref{eqn:Psi} that is suitable for reducing the numerical error from rounding, as exaggerated by dividing by a number close to zero.  We multiply by $R$ to get
\begin{equation}
R(\tau)\Psi^\prime(\tau) + \ell\sin\Psi(\tau) - \kappa\ell R(\tau) U(\tau) = 0, \label{eqn:Psialt}
\end{equation}
which we will use when $b$ is relatively small.  When we use \eqref{eqn:Psialt} in place of \eqref{eqn:Psi} we denote the changed $N$ by $N_1$ (and $F$ and $L$ similarly below) when specifying the difference is important, and we will generically refer to those objects without subscripts when no confusion is expected.

We will later approach this nonlinear problem with a Newton method, and we will need to use the Fr\'{e}chet derivative 
$$
F(\mathbf{v}) = \frac{d N}{d\mathbf{v}}(\mathbf{v}).
$$
Since $\mathbf{v}$ has several components, and some of those components involve derivatives with respect to $\tau$, we introduce the differential operator 
$$
D = \frac{d\, }{d\tau}, 
$$
which is applied in a block fashion to $\mathbf{v}$ so that $R^\prime(\tau) = [D\,\, 0\,\, 0\,\, 0]\mathbf{v}$, for example.  We will also have need to use an operator version of function evaluation.  We denote $D^0_\tau$ to be this operator, so that $D^0_1R = R(1)$.  With this in hand, we compute
\begin{equation} \label{eqn:N=0}
F(\mathbf{v} ) = 
\begin{bmatrix}
	D & 0 & \ell\sin\Psi & -\cos\Psi \\
	0 & D & -\ell\cos\Psi & -\sin\Psi \\
	\frac{-\ell\sin\Psi }{R^2} & -\kappa\ell & D + \frac{\ell\cos\Psi }{R} & \frac{\sin\Psi}{R} - \kappa U \\
	D^0_{-1} & 0 & 0 & 0 \\
	D^0_{1} & 0 & 0 & 0 \\
	0 & 0 & D^0_{-1} &  0 \\
	0 & 0 & D^0_{1} &  0 
\end{bmatrix}
\mathbf{v}.
\end{equation}
We will have need to solve linear systems based on the definition $F(\mathbf{v}) := L\mathbf{v}$.  For $F_1$ and $L_1$, we merely change the third row to 
\begin{equation}
\begin{bmatrix}
D\Psi - \kappa\ell U & -\kappa\ell R & RD + \ell\cos\Psi & \sin\Psi - \kappa UR
\end{bmatrix}.
\end{equation}

\subsection{ The problem {\bf P2} }
\label{subsecP2}

If we replace the symmetric and simply connected surfaces that are the image of the disk considered in Subsection~\ref{subsecP1} with doubly connected interfaces that are the image of an annulus, the scaling argument is preserved.  There are differences though, as the arclength $s=0$ does not correspond to the radius being zero.  We reuse \eqref{eqn:R}-\eqref{eqn:Psi} and \eqref{eqn:N=0}, but we have boundary conditions
\begin{eqnarray}
	R(1) - b &=& 0, \label{eqn:bcRb2} \\
	R(-1) - a &=& 0, \label{eqn:bcRa} \\
	\Psi(1) - \psi_b &=& 0, \label{eqn:bcPsib2} \\
	\Psi(-1) - \psi_a &=& 0, \label{eqn:bcPsia}
\end{eqnarray}
where we introduce $a\in(0,b)$ and $\psi_a\in[-\pi,\pi]$.  This changes the definition of $N$, but not that of $F$.  Notice that in this case there is no singularity as $R \ne 0$.  If $|\psi_a|, |\psi_b| \leq \pi/2$, then the resulting solution curves are graphs and we have an annular capillary tube, with the inner wall at the radius $a$ and the outer wall at the radius $b$ \cite{EKT2004}.  If this restriction of the inclination angles is removed, then the solution curve can pass through one or both of the ``walls''  \cite{Treinen2012} and is primarily useful as a component interface in a more complicated configuration, such as a floating drop problem \cite{ElcratTreinen2005} or a floating ball problem \cite{MccuanTreinen2013}.  In either case, these solution curves form generating curves with total arclength $2\ell$.  The global behavior can become quite complicated \cite{BagleyTreinen2018}, though we specify boundary conditions that avoid the self-intersections appearing in that work.

\subsection{ The problem {\bf P3} }
\label{subsecP3}

Our final prototype problem is the lower dimensional problem.  Here we do not have a radial symmetry, however the problem is quite similar to {\bf P2}.   Here we replace \eqref{drds}-\eqref{dpsids} with
\begin{eqnarray}
	\frac{dx}{ds} &=& \cos\psi, \label{dxds}\\
	\frac{du}{ds} &=& \sin\psi, \label{duds2}\\
	\frac{d\psi}{ds} &=& \kappa u, \label{dpsids2}
\end{eqnarray}
where $x$ represents a horizontal displacement.  There are two main ways to interpret this problem.  The first is that the solutions form a generating curve that is extended in the $y$ direction as a ruled surface.  One may take this extension in the $y$ direction as over a unit interval, or over any other interval, including the whole infinite interval.  In these cases the meridinal curvature is $\frac{d\psi}{ds}$ and the latitudinal curvature is zero.  The other interpretation is that this is quite simply a lower dimensional problem and there is only one curvature to consider.   For this reason we often call this the lower-dimensional problem, though it has other important interpretations.

This system is also known as Euler's elastica, as well as the nonlinear pendulum equation.  See \cite{Wente2006}, \cite{MR1368401}, and using Chebfun \cite{TrefethenBirkissonDriscoll2018}.  As in the previous case, the global behavior can become quite complicated \cite{Mccuan2022}, though we specify boundary conditions that avoid the self-intersections studied in that work.

As before, we scale the problem to find
\begin{eqnarray}
	X^\prime(\tau) - \ell\cos\Psi(\tau) &=& 0, \label{eqn:X}\\
	U^\prime(\tau) - \ell\sin\Psi(\tau) &=& 0, \label{eqn:U3}\\
	\Psi^\prime(\tau)  - \kappa\ell U(\tau) &=& 0. \label{eqn:Psi3}
\end{eqnarray}
where we define $X(\tau):= x(\ell\tau) = x(s)$ and the other functions are defined as above.   We use this, as before, to define $N(\mathbf{v}) = \mathbf{0}$ with $\mathbf{v} = [X\,\, U\,\, \Psi\,  \ell]^T$,  including  the boundary conditions
\begin{eqnarray}
	X(1) - b &=& 0, \label{eqn:bcXb} \\
	X(-1) - a &=& 0, \label{eqn:bcXa} \\
	\Psi(1) - \psi_b &=& 0, \label{eqn:bcPsib3} \\
	\Psi(-1) - \psi_a &=& 0, \label{eqn:bcPsia3}
\end{eqnarray}
with $-\infty < a < b<\infty$ and $-\pi \leq \psi_a, \psi_b \leq \pi$.  Then we compute
\begin{equation} \label{eqn:2D_N=0}
	F(\mathbf{v} ) = 
	\begin{bmatrix}
		D & 0 & \ell\sin\Psi & -\cos\Psi \\
		0 & D & -\ell\cos\Psi & -\sin\Psi \\
		0 & -\kappa\ell & D  &  - \kappa U \\
		D^0_{-1} & 0 & 0 & 0 \\
		D^0_{1} & 0 & 0 & 0 \\
		0 & 0 & D^0_{-1} &  0 \\
		0 & 0 & D^0_{1} &  0 
	\end{bmatrix}
	\mathbf{v}
\end{equation}
for this problem.  We again will use the notation $F(\mathbf{v}) = L\mathbf{v}$.

\section{The spectral method }
\label{Cheby}

The basis of the approach we take has its roots in Trefethen's monograph \cite{Trefethen2000}, though the treatment of boundary conditions has advanced since that work appeared.  The discretization of the differential operators is in the block form, as they appear above, and as is described in Driscoll and Hale \cite{DriscollHale2016} and then Aurentz and Trefethen \cite{AurentzTrefethen2017}.  Here we briefly summarize the ideas from those works.  Chebyshev differentiation matrices can be realized as the linear transformation between the data points corresponding to the interpolating polynomials for a function $f$ and its derivative $f^\prime$, where the data is sampled at Chebyshev grid points $x_j = \cos(\theta_j)\in[-1,1]$ where the angles $\theta_j$ are equally spaced angles over $[0,\pi]$.  Of course, these grid points do not need to be fixed, as multiple samplings can be used.  If the matrix is not square, then the number of rows being less than the number of columns can be seen as down-sampling.

We will treat the nonlinearity of $N$ with a Newton method.  Here we present the basic ideas with the backdrop of $\mathbf{P1}$.  We will postpone comments on initial guesses and other points of interest for the next section, and we will also discuss how to modify these basic ideas to the other problems in later sections.

The basic building blocks of  the nonlinear equation are based on $D^0$ and $D$, which we implement using the Chebfun commands
\begin{verbatim}
D0 = diffmat([n-1 n],0,X);
D1 = diffmat([n-1 n],1,X);
\end{verbatim}
where $X = [-1,1]$,  $n$ is the number of Chebyshev points we are using, and the input 0 or 1 indicates the number of derivatives.  Since $\texttt{D0}$  is rectangular, it becomes a $(n-1)\times n$ identity matrix interpreted as a dense ``spectral down-sampling'' matrix implemented as interpolating on an $n$-point grid followed by sampling on an $(n-1)$-point grid.    To sparsely build $N$ and $L$ we use the components

\begin{verbatim}
Z0 = sparse(n-1,n);
D01 = spalloc(n-1, 3*n + 1, n*(n-1));
D02 = D01, D03 = D01, D11 = D01, D12 = D01, D13 = D01;
D01(1:n-1, 1:n) = D0, D11(1:n-1, 1:n) = D1;
D02(1:n-1, n+1:2*n) = D0, D12(1:n-1, n+1:2*n) = D1;
D03(1:n-1, 2*n+1:3*n) = D0, D13(1:n-1, 2*n+1:3*n) = D1;
dT1n1 = sparse(1,3*n+1), dT1p1 = dT1n1;
dT3n1 = dT1n1, dT3p1 = dT1n1;
dT1n1(1) = 1, dT1p1(n) = 1, dT3n1(2*n+1) = 1, dT3p1(end-1) = 1;
\end{verbatim}
where we will construct our operator based on the building block of multiplying a linear operator times our vector $\mathbf{v}$ of solution components.  We see $\texttt{D01}$ allocates the appropriate sparse matrix, and then that portion of the matrix that corresponds to multiplying by the $R$ component of  $\mathbf{v}$ is set to $\texttt{D0}$ for function evaluation, and $\texttt{D11}$ uses $\texttt{D1}$ for the differentiation of $R$ there.  Then $\texttt{D02}$ evaluates $U$, $\texttt{D12}$ takes the derivative of $U$, and the same pattern is used for $\texttt{D03}$ and $\texttt{D13}$ with $\Psi$.  Then $\texttt{dT1n1, dT1p1, dT3n1}$, and $\texttt{dT3p1}$ evaluate $R$ and $\Psi$ at negative 1 and positive 1 respectively.

Then $N$ and $L$ are given by
\begin{verbatim}
N = @(v) [ D11*v - v(end).*cos(D03*v)
           D12*v - v(end).*sin(D03*v)
           D13*v + v(end).*sin(D03*v)./(D01*v) - kappa*v(end).*D02*v
           dT1n1*v + b; dT1p1*v - b
           dT3n1*v + psib; dT3p1*v - psib ];
\end{verbatim}
and
\begin{verbatim}
L = @(v) [ D1, Z0, spdiags(v(end)*sin(D03*v) ,0 ,n-1 ,n-1 )*D0, -cos(D03*v)
           Z0, D1, spdiags(-v(end)*cos(D03*v) ,0 ,n-1 ,n-1 )*D0, -sin(D03*v)
           spdiags(-v(end)*sin(D03*v)./((D01*v).^2) ,0 ,n-1 ,n-1 )*D0, ...
           -kappa*v(end)*D0, D1 + (spdiags(v(end)*cos(D03*v),0,n-1,n-1))*D0, ...
           sin(D03*v)./(D01*v) - kappa*D02*v
           dT1n1; dT1p1; dT3n1; dT3p1 ];
\end{verbatim}
where we use these version when $b > 1$.  This formulation works well for larger radii problems, as it does not exaggerate the numerical error for large $R$ values.

Here we take care that the number of Chebyshev points are chosen so that there is no evaluation of  $r = 0$.  This removable singularity is avoided in exactly the same way Trefethen described in \cite{Trefethen2000} when discussing the radial form of Laplace's equation.

If $b \leq 1$ then we use the formulations based on \eqref{eqn:Psialt} where we have no excessive numerical error due to rounding near the singularity.  Then $N_1$ and $L_1$ are given by
\begin{verbatim}
N1 = @(v) [ D11*v - v(end).*cos(D03*v)
            D12*v - v(end).*sin(D03*v)
            (D01*v).*(D13*v) + v(end).*sin(D03*v) - kappa*v(end).*(D02*v).*(D01*v)
            dT1n1*v + b; dT1p1*v - b
            dT3n1*v + psib; dT3p1*v - psib  ];
\end{verbatim}
and
\begin{verbatim}
L1 = @(v) [ D1, Z0, spdiags(v(end)*sin(D03*v),0,n-1,n-1)*D0, -cos(D03*v)
            Z0, D1, spdiags(-v(end)*cos(D03*v),0,n-1,n-1)*D0, -sin(D03*v)
            spdiags(D13*v - kappa*v(end).*(D02*v),0,n-1,n-1)*D0, ...
            spdiags(-kappa*v(end)*(D01*v),0,n-1,n-1)*D0, ...
            spdiags(D01*v,0,n-1,n-1)*D1 + spdiags(v(end)*cos(D03*v),0,n-1,n-1)*D0, ...
            sin(D03*v) - kappa*(D02*v).*(D01*v)
            dT1n1; dT1p1; dT3n1; dT3p1 ];
\end{verbatim}

The next steps are to construct initial guesses (see Section~\ref{Study1}) and use Newton's method.  To initialize this we pick a tolerance $\texttt{tol\_newton}$ and we use the relative error measured by the Frobeneous norm.  We found no practical difference when the infimum norm was used instead.   For all of the examples in this paper, we used $\texttt{tol\_newton = 1e-13}$.   The basic loop is
\begin{verbatim}
while res_newton > tol_newton
      dv = L(v)\N(v);
      v = v - dv;
      res_newton = norm(dv,'fro')/norm(v,'fro');
end
\end{verbatim}

If the Newton's method fails to converge within a specified number of iterations, we then increase the number of Chebyshev points by four and we sample the current state of the iteration onto the new Chebyshev points, reinitialize the operators $N$ and $L$ as above, and we enter another loop for Newton's method.  Further, we include this process in an outer  loop that also tests the relative error of the iterates even if the Newton's method converges.  We use the Frobenious norm to compute $\texttt{res\_bvp}$ as $||N(v)||/||v||$.  If this residual is greater than the prescribed tolerance, then we also increase the number of Chebyshev points by four as just described.  We used a prescribed tolerance of $\texttt{tol\_bvp = 1e-12}$ in all of the examples in this paper.  It may also be that there is excessive polynomial oscillation across the Chebyshev points, and there are not enough of these points to accurately resolve the solution.  We will explain this further below, but we note here that in this case then we increase the number of Chebyshev points by $2(n - 1) - 1$ instead.  It may be that for some problems this adaptive parameter of the number of new points to add benefits from a customized choice.   In the sense just described, our algorithm is adaptive.  

If the boundary conditions are specified so that the solution curve will not be a graph over a base-domain, such as when $\psi_b > \pi/2$, then we use a method of continuation.  We harvest the sign of $\psi_b$ and we solve ten problems with linearly spaced angles between $\pi/2$ and $|\psi_b|$ using the previously converged solution at each step as an initial guess for the next step.  We preserve the number of needed Chebyshev points from the adaptive algorithm while moving from one step to the next.  If $\psi_b < 0$, then we reflect about the $r$-axis.  In problems $\mathbf{P2}$ and $\mathbf{P3}$, if needed, we do this for the inclination angles at $a$ and $b$ simultaneously.

\begin{figure}[hb]
	\centering
	\scalebox{0.35}{\includegraphics{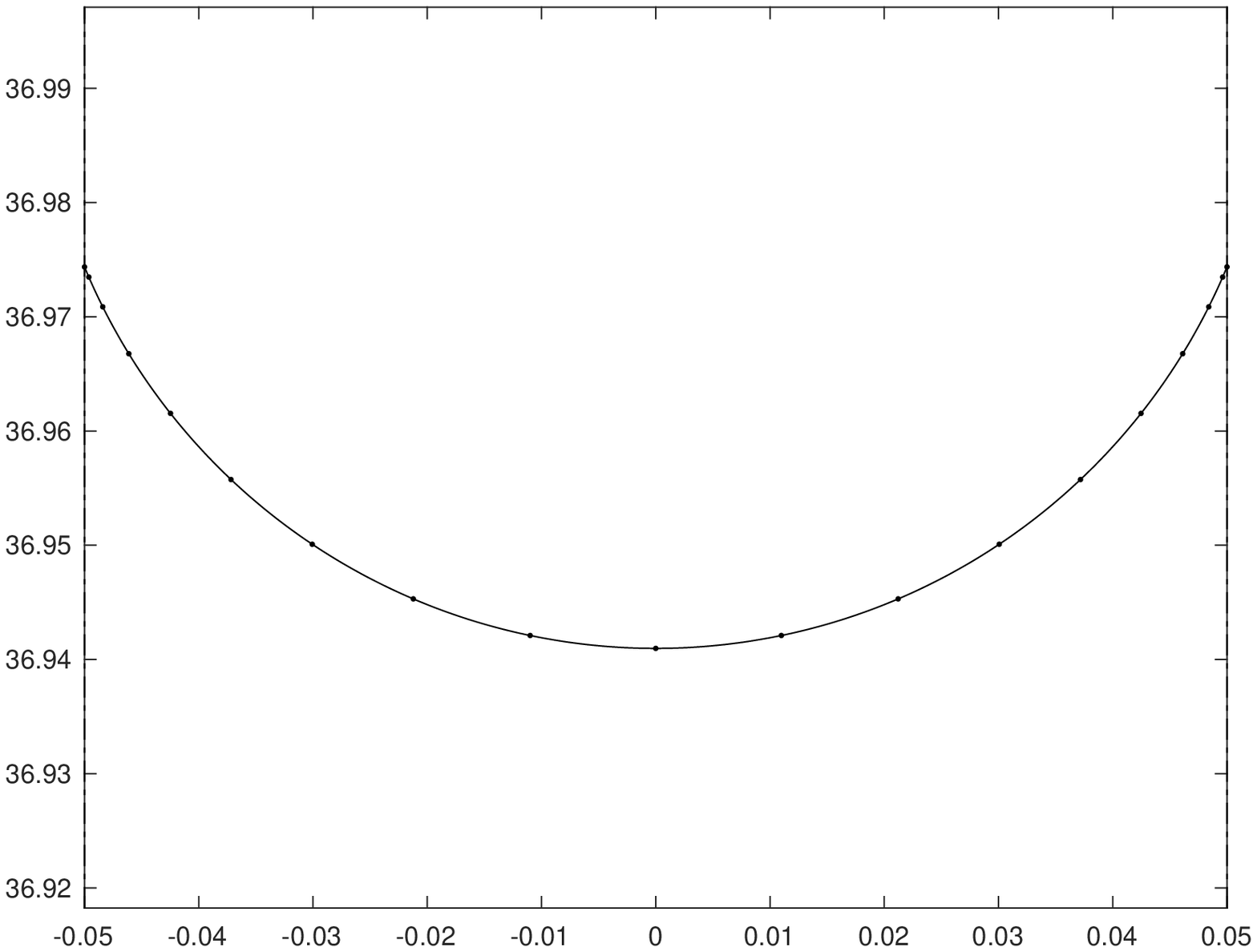}}
	\scalebox{0.35}{\includegraphics{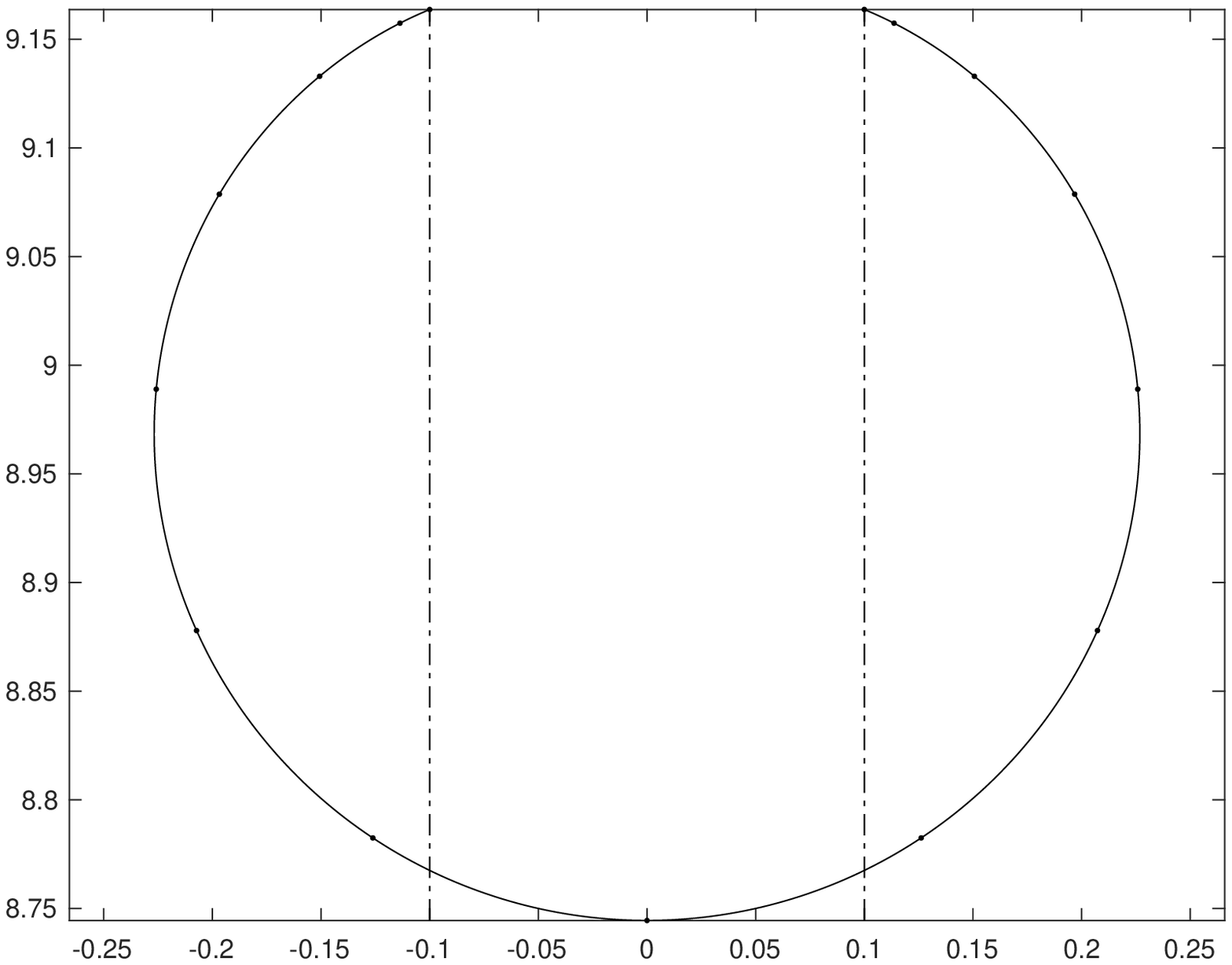}}
	\caption{A capillary tube with radius $b = 0.05$ and inclination angle $\psi_b = 3\pi/8$ (L) and with radius $b = 0.01$ and  inclination angle $\psi_b = 7\pi/8$ (R).  In all of our figures we include vertical lines to indicate the radius (or radii) of interest.  In some figures this is obscured by the bounding boxes of the figure.}
	\label{fig:cap_small}
\end{figure}

\section{{\bf P1} Study}
\label{Study1}

In order  to start Newton's method, we need an initial guess for the solution.    With this in mind, we begin our study of symmetric capillary surfaces which are the image of a disk, and we first consider some details known about the surfaces themselves.    For the moment we will restrict our attention to $|\psi_b|\leq \pi/2$ as in Figure~\ref{fig:cap_small}~(L), and if  $|\psi_b| >  \pi/2$ as in Figure~\ref{fig:cap_small}~(R) we will use the method of continuation as described above.  Further, we take $\psi_b \geq 0$, and reflect if needed.

\begin{table}[hb]
	\centering
	\begin{tabular}{||l|l|l|l||}
		\cline{1-4}
		$b$ & $\psi_b$ & $n$ & Elapsed time \\
		\hline \hline
		1 & $\pi/6$ & 15 & 0.0152 \\
		\hline
		1 & $\pi/3$ & 15 & 0.0134  \\
		\hline
		1 & $\pi/2$ & 15 & 0.0050 \\
		\hline
		1 & $2\pi/3$ & 15 & 0.0358  \\
		\hline
		1 & $\pi$ & 15 & 0.0284\\
		\hline
	\end{tabular}
	\caption{Solutions that are the image of a disk with radius $b = 1$, and a range of inclination angles, where the number of Chebyshev points $n$ and the elapsed time in seconds is shown.}
	\label{table_disc_1}
\end{table}

\begin{table}[t]
	\centering
	\begin{tabular}{||c|l|>{\columncolor{gray!15}}l|l|>{\columncolor{gray!15}}l||l|>{\columncolor{gray!15}}l||}
		\hline
		$b$ &  \multicolumn{2}{|c|}{$\psi_b$} & \multicolumn{2}{|c|}{$n$} & \multicolumn{2}{|c||}{Elapsed time }  \\
		\hline \hline
		0.05 &    $3\pi/8$ &  $7\pi/8$ & 15  & 15  & 0.0603   &  0.0422\\
		\hline
		0.10 & $3\pi/8$ &  $7\pi/8$ & 15   & 15  &  0.0169   &  0.0367\\
		\hline
		0.50 &   $3\pi/8$ &  $7\pi/8$ & 15   & 15  & 0.0054   & 0.0376\\
		\hline
		1 &    $3\pi/8$ &  $7\pi/8$ & 15   & 15  &  0.0056  &  0.0362\\
		\hline
		2 &    $3\pi/8$ &  $7\pi/8$ & 15   & 15  &  0.0398  &  0.1043\\
		\hline
		10 &   $3\pi/8$ &  $7\pi/8$ & 27   & 75  & 0.0452 &  0.4510\\
		\hline
		20 &   $3\pi/8$ &  $7\pi/8$ & 75   & 127  &  0.1613 & 1.2785\\
		\hline
	\end{tabular}
	\caption{Solutions that are the image of a disk with a range of radii $b$, and two choices of inclination angles shown in sub-columns. Also shown are the number of Chebyshev points $n$ and the elapsed time in seconds.  The shaded sub-columns indicate the corresponding data. }
	\label{table_disc_2}
\end{table}

The asymptotic limit of the solution as $b\rightarrow 0$, as given by the generating curve, is known to be a circular arc.  Laplace gave (without proof) an estimate of the height of the interface on the axis of symmetry \cite{MR0265115}, and Siegel established this rigorously in 1980 \cite{Siegel1980}.  See Finn's monograph \cite{ecs} for the details of what we will use here.  Finn describes two circular arcs, with one strictly below the generating curve away from the $u$-axis, and another strictly above.  We have found that the circular arc above is easy to compute from the boundary conditions, and up to vertical translations to adjust for estimates of the height of the tip, fairly robust.  If the container has a small radius, the more precise formulation of the average radius of the two circular arcs is used with the vertical translation given by Laplace's formula.  In all cases we will use the notation of $u_0$ for the height of the tip.  Laplace's formula gives
\begin{eqnarray}
	u_0 &\approx& \frac{2\cos\gamma}{\kappa b} - \frac{b}{\cos\gamma} + \frac{2b}{3}\frac{1 - \sin^3\gamma}{\cos^3\gamma}  \nonumber\\
	&=& \frac{2\cos\gamma}{\kappa b} - \frac{b\cos\gamma}{3}\frac{ 1 + 2\sin\gamma}{(1 + \sin\gamma)^2},
\end{eqnarray}
where the contact angle  is $\gamma = \pi/2 - \psi_b$.
Finn shows that the two sectional curvatures are equal to $\kappa b^2 u_0/2$, and so we define the lower circular arc $\Sigma^{(1)}$ and the corresponding function $u^{(1)}(r)$, centered on the $u$-axis with $u^{(1)}(0)$ to be the height of the  tip and radius $R_1 = 2/(\kappa b^2 u_0)$.  He finds this circular arc to be below the solution curve.  Then he defines $\Sigma^{(2)}$ and the corresponding function $u^{(2)}(r)$, centered on the $u$-axis with $u^{(2)}(0)$ to be the height of the  tip, and  so that the inclination angle matches the prescribed angle at $b$.  He establishes the fact that this curve is above the solution curve.  We then have formulas
$$
u^{(1)}(r) = u_0 + R_1 - \sqrt{R_1^2 - r^2},
$$
and
$$
u^{(2)}(r) = u_0 + R_2 - \sqrt{R_2^2 - r^2},
$$
where we derive that $R_2 = b\csc\psi_b$.

If $b < 1$, we use the averages of these two radii to construct a circular arc that is midway between the upper and lower estimates for use in our initial guess for the height $U$ as a function of $R$, and we use the initial height $u_0$ as above.  In the cases when $b \geq 1$ we have found $u^{(2)}(r)$ with $u_0 = 2/R_2$  to be a better estimate for our initial guess.    To build this we of course need estimates for the other quantities, and we always assume the initial guess to be a circular arc.  We then have $\Psi_0(\tau) = \tau\psi_b$ as a function over the domain of Chebyshev points, and $R(\tau) = \hat R\sin(\Psi_0(\tau))$  where $\hat R$ is whichever radius we are using for that case.  Finally, we use $\ell_0 = \abs{\hat R\psi_b}$ for the arclength.

In all of these cases, to account for an angle $\psi_b < 0$,  we multiply the relevant corresponding initial guess components by $sign(\psi_b)$.

It can be that the iterates in Newton's method stray away from physically relevant solutions, and sometimes get stuck there.  To prevent this from happening two barriers are used in the iteration process.  First, inside of the Newton's method loop, if $|\Psi | > 3.5$, then those corresponding entries of the computational vector are reset to $\pm\pi$.  Physically, we never have angles greater than $\pi$, so this restriction removes these problematic iterations.   Since the solution curves have at most one inflection point, and none in the case we are considering, we measure the oscillations of $\Psi$ and if  
\begin{equation}\label{eqn:oscillation}
\texttt{length(find(diff(sign(diff((S2*v(2*n+1:3*n))))))) >= 2}
\end{equation}
then for that iterate, we linearly interpolate $\Psi$ between the boundary conditions.  In this case, that amounts to linearly interpolating from $-\psi_b$ to $\psi_b$.  

The other case we can encounter is that there are not enough Chebyshev points to accurately resolve the solution.  This typically appears in the outer  loop, outside of the Newton's method loop, and it shows up as excessive numerical oscillation of the approximating polynomial.  So there in that outer loop we again test \eqref{eqn:oscillation}, and if this holds, then we simply increase the number of Chebyshev points by $2(n-1) -1$, as described above.

We include some results from running this code in Tables~\ref{table_disc_1}-\ref{table_disc_2}.  The number of Chebyshev nodes is reasonably small, so the matrices are of size $(3n + 1)\times(3n + 1)$, and the code runs quite quickly, even when a rather extreme radius is chosen.

\section{{\bf P2} Study}
\label{Study2}

In our study of  bounded  simply connected symmetric capillary surfaces, we begin by using an initial guess based on the ideas in the last section.  

If $\psi_a\psi_b <0$, then the solution will be ``u-shaped''.  The upward oriented case when $\psi_b > 0 $ is reflected about the radial axis when $\psi_b < 0$.  Then we center a circular approximation at $(b - a)/2$ and we derive the radius to be 
$$
\hat R = \frac{b - a}{\sin(-\psi_a) + \sin(\psi_b)},
$$
so that it meets $a$ and $b$ at angles $\psi_a$ and $\psi_b$ respectively.
Our arclength approximation is  $\ell_0  = \hat R|\psi_b - \psi_a|$,  we use $u_0 = 1/\hat R$, and we linearly interpolate the inclination angle by 
$$
\Psi_0(\tau) = \frac{1 - \tau}{2}\psi_a + \frac{1 + \tau}{2}\psi_b.
$$
The circular approximation of the solution curve is then given by
$$
R_0(\tau) = \frac{a+b}{2} + \hat R \sin(\Psi_0(\tau)),
$$
and
$$
U_0(\tau) = u_0 + \hat R(1 - \cos(\Psi_0(\tau))).
$$

If $\psi_a\psi_b \geq 0$, then we have an ``s-shaped'' curve.  Again, the case with upward oriented right half of the curve is given by $\psi_b > 0$, and if $\psi_b < 0$ we simply reflect the construction   about the radial axis.  Our approximation of the interface is not strictly circular, but shares some of the same concepts.  We set
$$
\hat R = \frac{b - a}{\sin(\psi_a) + \sin(\psi_b)},
$$
with $\ell_0  = \hat R(\psi_b + \psi_a)$, then
$$
\Psi_0(\tau) = \frac{1 - \tau}{2}\psi_a + \frac{1 + \tau}{2}\psi_b,
$$
and
$$
R_0(\tau) = \frac{a+b}{2} + \hat R \sin(\Psi_0(\tau)),
$$
are virtually the same.  The most substantial difference is that we define 
$$
U_0(\tau) = \frac{\tau}{b-a},
$$
which changes sign and is linear.  While this is not a very accurate approximation, in practice it has shown to be fairly robust as an initial guess in Newton's method.

We modify the above formulations of $N$ and $L$ as
\begin{verbatim}
N = @(v) [ D11*v - v(end).*cos(D03*v)
           D12*v - v(end).*sin(D03*v)
           D13*v + v(end).*sin(D03*v)./(D01*v) - kappa*v(end).*D02*v
           dT1n1*v - a; dT1p1*v - b
           dT3n1*v - psia; dT3p1*v - psib ];
\end{verbatim}
and
\begin{verbatim}
L = @(v) [ D1, Z0, spdiags(v(end)*sin(D03*v),0,n-1,n-1)*D0, -cos(D03*v)
           Z0, D1, spdiags(-v(end)*cos(D03*v),0,n-1,n-1)*D0, -sin(D03*v)
           spdiags(-v(end)*sin(D03*v)./((D01*v).^2),0,n-1,n-1)*D0, ...
           -kappa*v(end)*D0, D1 + (spdiags(v(end)*cos(D03*v),0,n-1,n-1))*D0, ...
           sin(D03*v)./(D01*v)  - kappa*D02*v
           dT1n1; dT1p1; dT3n1; dT3p1 ];
\end{verbatim}
and
\begin{verbatim}
N1 = @(v) [ D11*v - v(end).*cos(D03*v)
            D12*v - v(end).*sin(D03*v)
            (D01*v).*(D13*v) + v(end).*sin(D03*v) - kappa*v(end).*(D02*v).*(D01*v)
            dT1n1*v - a; dT1p1*v - b
            dT3n1*v - psia; dT3p1*v - psib  ];
\end{verbatim}
and
\begin{verbatim}
L1 = @(v) [ D1, Z0, spdiags(v(end)*sin(D03*v),0,n-1,n-1)*D0, -cos(D03*v)
            Z0, D1, spdiags(-v(end)*cos(D03*v),0,n-1,n-1)*D0, -sin(D03*v)
            spdiags(D13*v - kappa*v(end).*(D02*v),0,n-1,n-1)*D0, ...
            spdiags(-kappa*v(end)*(D01*v),0,n-1,n-1)*D0, ...
            spdiags(D01*v,0,n-1,n-1)*D1 + spdiags(v(end)*cos(D03*v),0,n-1,n-1)*D0, ...
            sin(D03*v) - kappa*(D02*v).*(D01*v)
            dT1n1; dT1p1; dT3n1; dT3p1 ];
\end{verbatim}
and the rest of the algorithm proceeds without modification. 

\begin{figure}[!h]
	\centering
	\scalebox{0.35}{\includegraphics{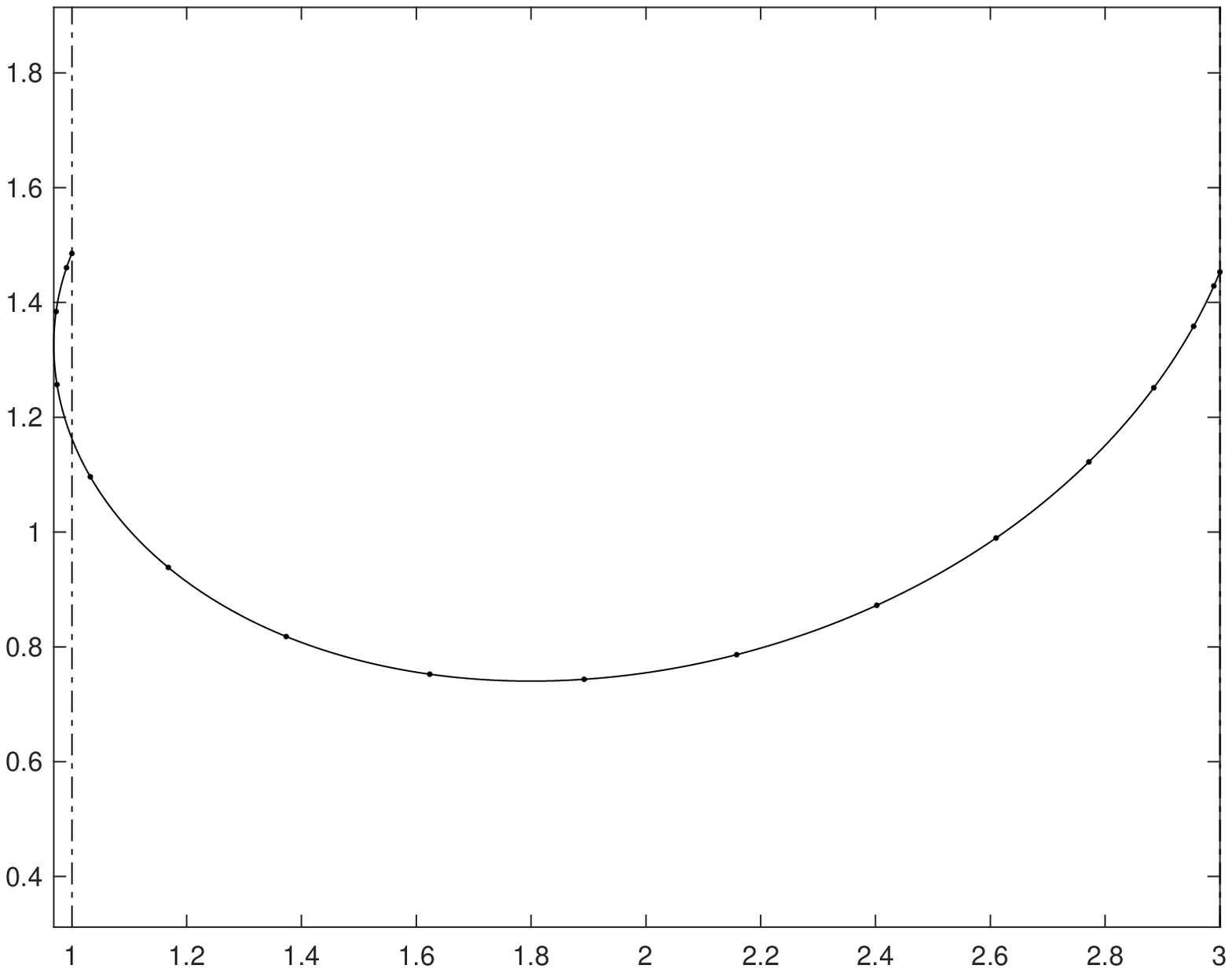}}
		\scalebox{0.35}{\includegraphics{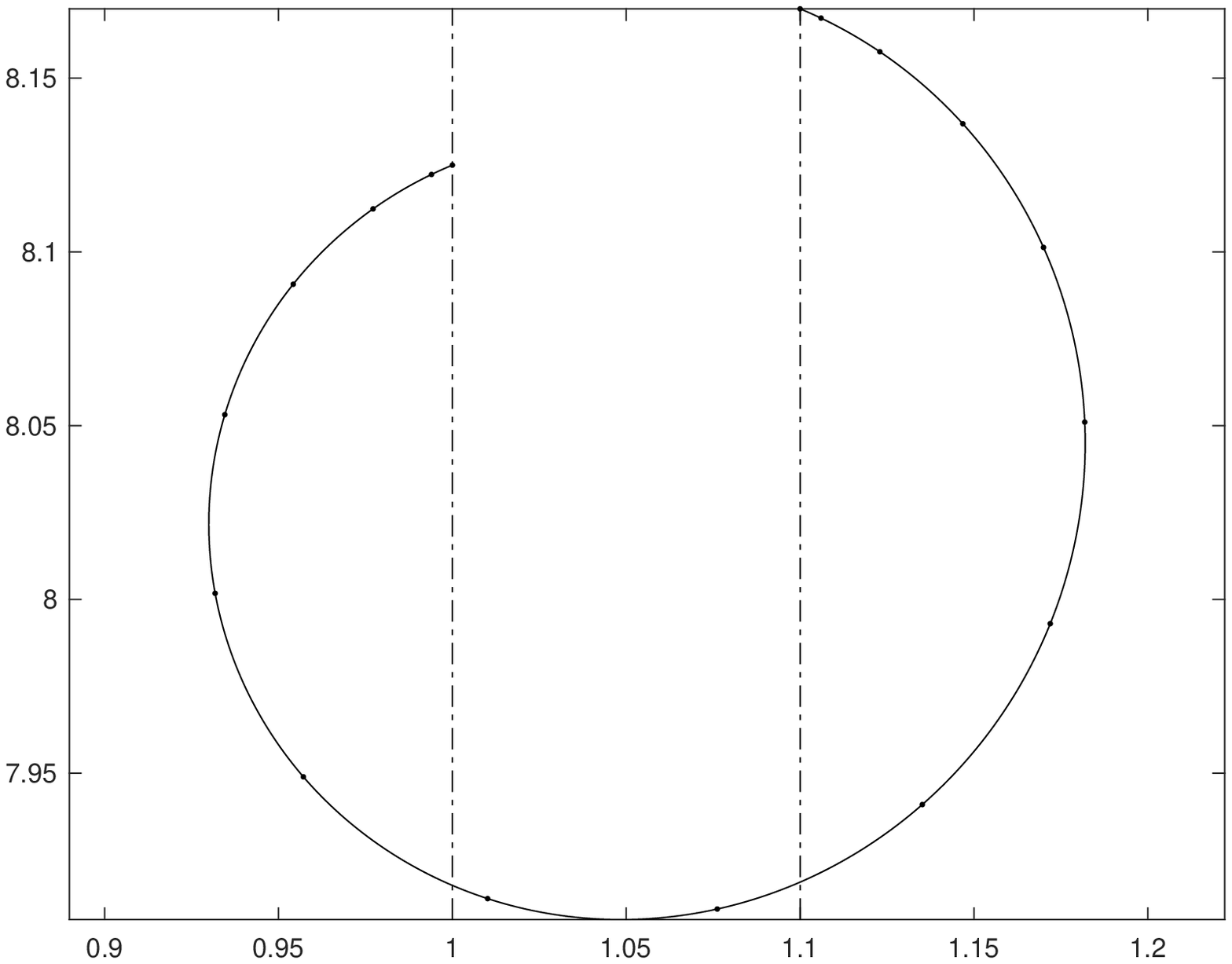}}
	\caption{A annular surface  with inner radius $a = 1$ and outer radius $b = 3$ and  inclination angles $\psi_a = -5\pi/8$ and $\psi_b = 3\pi/8$ there (L) and an annular surface  with inner radius $a = 1$ and outer radius $b = 1.1$ and  inclination angles $\psi_a = -7\pi/8$ and $\psi_b= 7\pi/8$ (R).}
	\label{fig:annular_1_3}
\end{figure}

In Figure~\ref{fig:annular_1_3} (L) we show a generating curve that has a vertical point on the left side of its span, but no vertical point on the right.  This type of configuration is typical for multi-component problems where there is an object or a different fluid ``floating'' in the center of a cylindrical container of radius $b$ and the contact of the ``exterior'' interface with that object or fluid is at radius $a$.  The physical (or mathematical) considerations may both determine $a$ and find $|\psi_a| > \pi/2$ in some cases.

In Figure~\ref{fig:annular_1_3} (R) we see an example where there are two vertical points which highlights the differences in the mathematical limits of interfaces as the radii get small.  Here, as the gap between $a$ and $b$ gets small, the differences in the curvature on the right and the left are apparent when compared to problem {\bf P1} when $b \rightarrow 0$ as in Figure~\ref{fig:cap_small}.

\begin{figure}[!h]
	\centering
	\scalebox{0.35}{\includegraphics{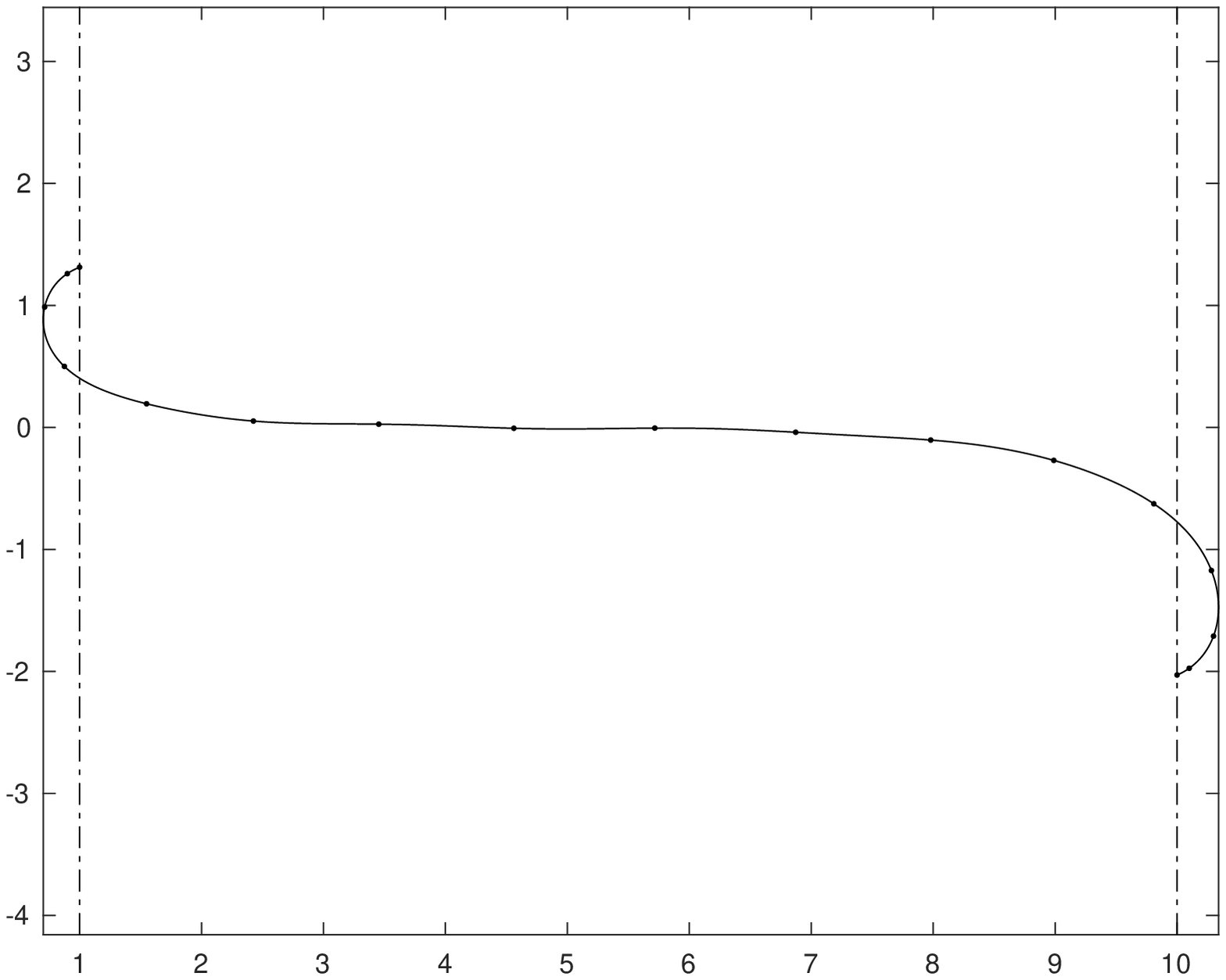}}
		\scalebox{0.35}{\includegraphics{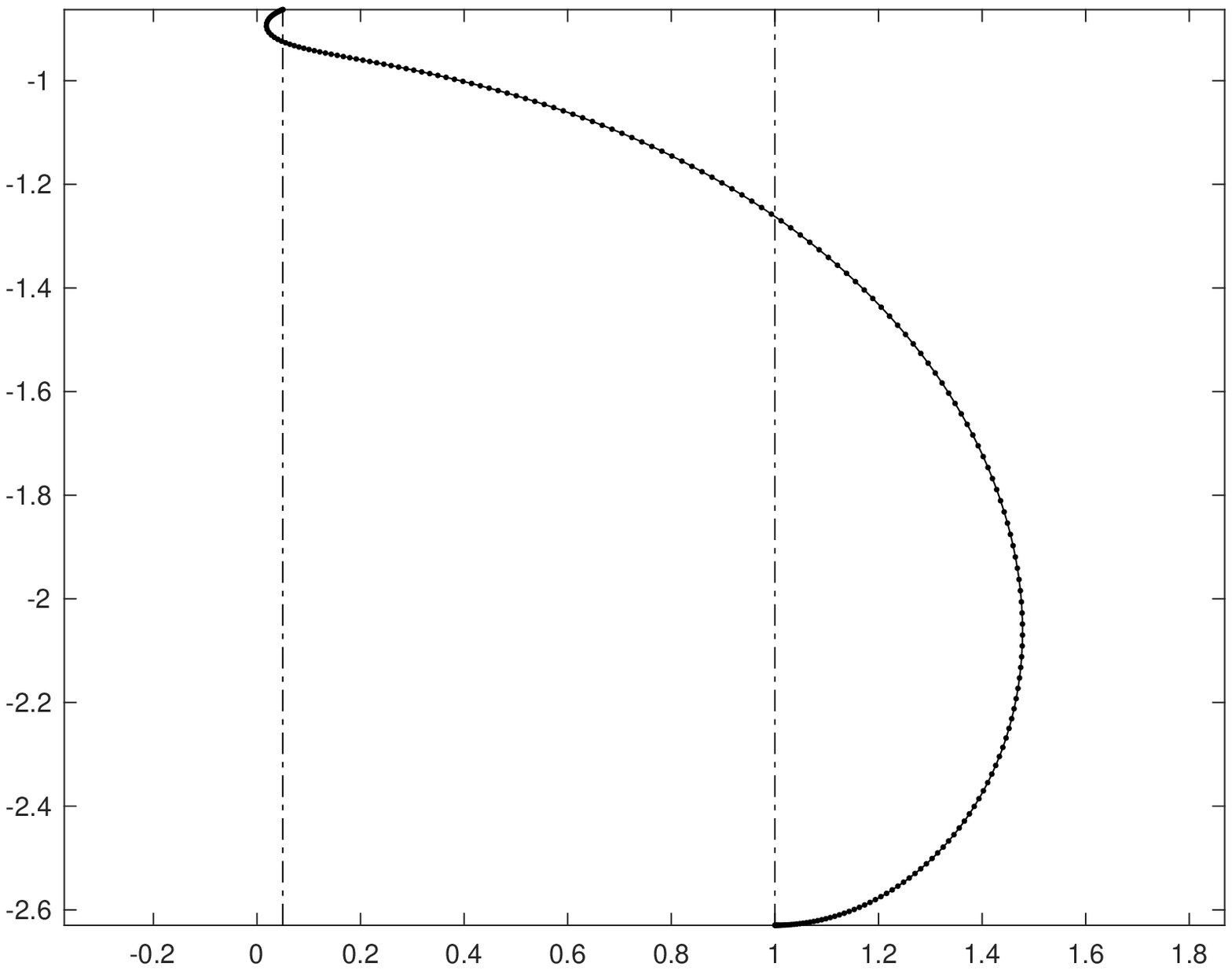}}
	\caption{A annular surface  with inner radius $a = 1$ and outer radius $b = 10$ and  inclination angles $\psi_a = -7\pi/8$ and $\psi_b = -7\pi/8$ there (L) and an annular surface  with radii $a = 0.05$ and  $b = 1$,  inclination angles $\psi_a = -7\pi/8$ and $\psi_b= \pi$ that was computed in 0.2242 seconds and finished with 195 Chebyshev points (R).}
	\label{fig:annular_1_10}
\end{figure}

\begin{figure}[h]
	\centering
	\scalebox{0.35}{\includegraphics{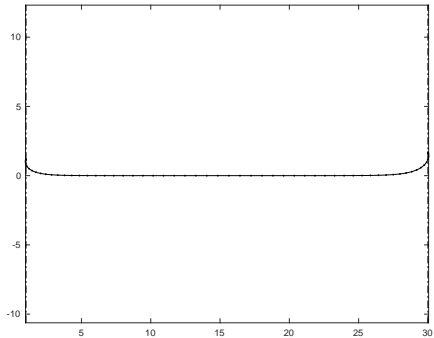}}
	\caption{A annular surface  with inner radius $a = 1$ and outer radius $b = 30$ and  angles $\psi_a = -5\pi/8$ and $\psi_b = 5\pi/8$ using 61 data points along the interface.}
	\label{fig:annular_1_30}
\end{figure}

In Figure~\ref{fig:annular_1_10} are examples of  ``s-shaped'' interfaces, also containing two vertical points.  On the left, we highlight a larger gap between $a$ and $b$ that used 17 Chebyshev points to resolve the interface to our requested precision.  This is then to be contrasted with Figure~\ref{fig:annular_1_30} where the gap between $a$ and $b$ is larger, but we needed 61 Chebyshev points to achieve our requested precision.  This last figure also gives a concrete example of where the algorithm presented here begins to have real advantages over the shooting method used previously  (say, in \cite{ElcratTreinen2005}).  While we will describe this advantage in more detail in the next section, the basic idea of this ``legacy'' algorithm is to use a point $(m,u_m)$ corresponding to the horizontal point in the middle of the region between $a$ and $b$.  (Clearly this only works for u-shaped interfaces, and the basic ideas have been modified to other configurations.)  Then one uses the ODEs \eqref{drds}-\eqref{dpsids} to propagate the curve to the left and right, stopping at the prescribed inclination angle there.  Then a nonlinear solver is used to match the radii that are produced by the construction just sketched with the prescribed radii $a$ and $b$ by varying the starting point $(m,u_m)$.  This shooting method  does not perform very well when the angle $\psi$ and $du/ds$ are both near zero for a large range of radii, meaning that there has been a barrier on the size of the gap in radii numerically studied up to now.  This also applies to {\bf P1} for large $b$.  The literature does not contain any discussion of this limitation that this author is aware of. 

On the right of Figure~\ref{fig:annular_1_10} we exhibit the beginnings of a multi-scale behavior.  Here $a = 0.05$  is chosen quite small, and $b = 1$ is moderate, with inclination angles $\psi_a = -7\pi/8$ and $\psi_b= \pi$.  This configuration converges fairly quickly, but it takes quite a few Chebyshev points.  If $\psi_a = -\pi$ is chosen, the code does not converge with any reasonable speed.  In general this phenomena does sometimes cause our code to fail, and further refinement is needed to treat these cases with such multi-scale features.

We include some results from running this code in Tables~\ref{table_ann_1}-\ref{table_ann_3}.  We include in these numerical studies results for the annular-type problems as well as the results for the lower-dimensional problem with the same data.  The number of Chebyshev nodes is reasonably small, so the matrices are of size $(3n + 1)\times(3n + 1)$, and the code runs quite quickly.

\begin{table}[h!]
	\centering
	\begin{tabular}{||l|c|c|l| >{\columncolor{gray!15}}l|l|>{\columncolor{gray!15}}l|l|>{\columncolor{gray!15}}l|l|>{\columncolor{gray!15}}l|l|>{\columncolor{gray!15}}l ||}
		\hline
		     \multicolumn{5}{||c|}{}& \multicolumn{4}{|c|}{Radially symmetric} & \multicolumn{4}{|c||}{Lower dimensional} \\
\hline
		$a$ & $b$ & $\psi_a$ & \multicolumn{2}{|c|}{$\psi_b$} & \multicolumn{2}{|c|}{$n$} & \multicolumn{2}{|c|}{Elapsed time} & \multicolumn{2}{|c|}{$n$} &  \multicolumn{2}{|c||}{Elapsed time} \\
		\hline \hline
     1 &   3 &         0    & $\pi/6$  & $-\pi/6$     &15&15&    0.0054  &  0.0062  &15&15&   0.0020  &  0.0029 \\
     1 &   3 &  $-\pi/6$     & $\pi/6$  & $-\pi/6$    &15&15&    0.0049  &  0.0048  &15&15&    0.0021  &  0.0026\\
     1 &   3 &   $-\pi/3$     & $\pi/6$   & $-\pi/6$ &15&15&    0.0050  &  0.0054   &15&15&  0.0023  &  0.0028 \\
     1 &   3 &     $-\pi/2$     & $\pi/6$   & $-\pi/6$ &15&15&    0.0061  &  0.0066   &15&15&    0.0035  &  0.0025 \\
     1 &   3 &     $-2\pi/3$     & $\pi/6$  & $-\pi/6$  &15&15&    0.0569  &  0.0524   &15&15&     0.0178  &  0.0191 \\
     1 &   3 &     $-\pi$     & $\pi/6$  & $-\pi/6$  &15&15&   0.0418  &  0.0361  &15&15&     0.0157  &  0.0157\\
          \hline
     1 &   3 &         0    &  $\pi/3$  & $-\pi/3$  &15&15&   0.0037  &  0.0046    &15&15&     0.0019  &  0.0023\\
     1 &   3 &   $-\pi/6$     &  $\pi/3$  & $-\pi/3$   &15&15&    0.0040  &  0.0040  &15&15&     0.0018  &  0.0018 \\
     1 &   3 &    $-\pi/3$     &  $\pi/3$   & $-\pi/3$  &15&15&      0.0039 &   0.0036   &15&15&    0.0015  &  0.0016\\
     1 &   3 &     $-\pi/2$     &  $\pi/3$   & $-\pi/3$  &15&15&    0.0037  &  0.0034    &15&15&     0.0016  &  0.0018\\
     1 &   3 &     $-2\pi/3$     &  $\pi/3$  &  $-\pi/3$    &15&15&     0.0351  &  0.0352   &15&15&   0.0180  &  0.0151\\
     1 &   3 &     $-\pi$     &  $\pi/3$   & $-\pi/3$  &15&15&    0.0398  &  0.0369   &15&15&    0.0202  &  0.0187\\
     \hline
     1 &   3 &         0      & $\pi/2$   & $-\pi/2$  &15&15&  0.0067   & 0.0054    &15&15&    0.0021  &  0.0025\\
     1 &   3 &   $-\pi/6$  &     $\pi/2$  & $-\pi/2$  &15&15&     0.0047  &  0.0046    &15&15&     0.0021  &  0.0021\\
     1 &   3 &    $-\pi/3$  &   $\pi/2$  & $-\pi/2$  &15&15&     0.0043  &  0.0044 &15&15&     0.0020  &  0.0021 \\
     1 &   3 &     $-\pi/2$    &   $\pi/2$  & $-\pi/2$  &15&15&     0.0042  &  0.0042    &15&15&      0.0019  &  0.0020 \\
     1 &   3 &     $-2\pi/3$   &    $\pi/2$ & $-\pi/2$   &15&15&   0.0373  &  0.0371  &15&15&   0.0172  &  0.0183\\
     1 &   3 &     $-\pi$   &    $\pi/2$ & $-\pi/2$   &15&15&     0.0393  &  0.0377  &15&15&     0.0152  &  0.0151\\
          \hline
     1 &   3 &         0    &  $2\pi/3$  &  $-2\pi/3$   &15&15&   0.0364  &  0.0411  &15&15&   0.0152  &  0.0154\\
     1 &   3 &   $-\pi/6$  &     $2\pi/3$  &  $-2\pi/3$  &15& 17&    0.0341  &  0.0339   &15 &15&    0.0148  &  0.0150 \\
     1 &   3 &    $-\pi/3$     &  $2\pi/3$  &  $-2\pi/3$  &15&15&   0.0342  &  0.0363   &15&15&    0.0199  &  0.0196\\
     1 &   3 &     $-\pi/2$    &   $2\pi/3$  &  $-2\pi/3$  &15&15&    0.0402  &  0.0440   &15&15&  0.0165 &   0.0171\\
     1 &   3 &     $-2\pi/3$   &    $2\pi/3$  &  $-2\pi/3$  &15&15&     0.0412  &  0.0442   &15 &15&     0.0163  &  0.0150\\
     1 &   3 &     $-\pi$    &   $2\pi/3$  &  $-2\pi/3$  &15&15&    0.0406  &  0.0400   &15 &15&   0.0219  &  0.0200 \\
          \hline
     1 &   3 &         0  &    $\pi$   &    $-\pi$  &15&15&    0.0450  &  0.0415    &15&15&     0.0172  &  0.0206 \\
     1 &   3 &    $-\pi/6$  &     $\pi$   &    $-\pi$ &15&15&     0.0432  &  0.0534     &15&15&    0.0286  &  0.0256\\
     1 &   3 &     $-\pi/3$   &    $\pi$   &    $-\pi$ &15&15&     0.0743  &  0.0547   &15&15&       0.0195  &  0.0208\\
     1 &   3 &     $-\pi/2$    &   $\pi$   &    $-\pi$ &15&15&   0.0457  &  0.0441   &15&15&   0.0165  &  0.0179\\
     1 &   3 &     $-2\pi/3$   &    $\pi$   &    $-\pi$ &15&15&    0.0494  &  0.0446   &15&15&    0.0185  &  0.0179\\
     1 &   3 &     $-\pi$  &     $\pi$   &    $-\pi$ &15&15&   0.0354  &  0.0396  &15&15&     0.0152  &  0.0175  \\
\hline
\end{tabular}
\caption{Solutions that are the image of an annulus and the lower-dimensional problem with radii $a$ and $b$, and a range of inclination angles shown in sub-columns. Also shown are the number of Chebyshev points $n$ and the elapsed time in seconds.  The shaded sub-columns indicate the corresponding data.}
\label{table_ann_1}
\end{table}

\begin{table}[h!]
	\centering
	\begin{tabular}{||l|c|c|l| >{\columncolor{gray!15}}l|l|>{\columncolor{gray!15}}l|l|>{\columncolor{gray!15}}l|l|>{\columncolor{gray!15}}l|l|>{\columncolor{gray!15}}l ||}
		\hline
		\multicolumn{5}{||c|}{}& \multicolumn{4}{|c|}{Radially symmetric} & \multicolumn{4}{|c||}{Lower dimensional} \\
		\hline
		$a$ & $b$ & $\psi_a$ & \multicolumn{2}{|c|}{$\psi_b$} & \multicolumn{2}{|c|}{$n$} & \multicolumn{2}{|c|}{Elapsed time} & \multicolumn{2}{|c|}{$n$} &  \multicolumn{2}{|c||}{Elapsed time} \\
		\hline \hline
 1 &   1.05 &    $-3\pi/8$  &     $3\pi/8$  &    $-3\pi/8$      & 15 & 15 &   0.0032  &  0.0029   & 15 & 15 &    0.0017   & 0.0022\\
1 &    1.10 &   $-3\pi/8$  &     $3\pi/8$  &    $-3\pi/8$      & 15 &  15 &  0.0049   & 0.0038    & 15 &  15 &   0.0023  &  0.0024 \\
1 &    1.50 &    $-3\pi/8$  &     $3\pi/8$  &    $-3\pi/8$      & 15 &  15 &  0.0041  &  0.0049    & 15 &  15 &   0.0026  &  0.0025\\
1 &    2.00 &    $-3\pi/8$  &     $3\pi/8$  &    $-3\pi/8$      & 15 &   15 &  0.0048  &  0.0045    & 15 &  15 &   0.0022  &  0.0024 \\ 
1 &   10.00 &    $-3\pi/8$  &     $3\pi/8$  &    $-3\pi/8$      &15& 15&   0.0083  &  0.0070    &15& 15&   0.0023  &  0.0025 \\
1 &   15.00 &     $-3\pi/8$  &     $3\pi/8$  &    $-3\pi/8$      & 23 &  15 &  0.0161  &  0.0065   &15& 15&    0.0021  &  0.0024\\
\hline
1 &    1.05 &    $-3\pi/8$  &     $7\pi/8$  &    $-7\pi/8$     &15& 15&  0.0311  &  0.0328   &15& 15&   0.0184  &  0.0167\\
1 &    1.10 &    $-3\pi/8$  &     $7\pi/8$  &    $-7\pi/8$      &15& 15&   0.0272   & 0.0264    &15&  15&   0.0158  &  0.0180 \\
1 &    1.50 &    $-3\pi/8$  &     $7\pi/8$  &    $-7\pi/8$      &15& 15&   0.0288  &  0.0310    &15& 15&   0.0185  &  0.0154 \\
1 &    2.00 &    $-3\pi/8$  &     $7\pi/8$  &    $-7\pi/8$      &15& 15&  0.0345  &  0.0318   &15& 15&   0.0166   & 0.0172 \\
1 &   10.00 &    $-3\pi/8$  &     $7\pi/8$  &    $-7\pi/8$      &15& 15&  0.0473 &   0.0470    &15& 15&   0.0167  &  0.0155 \\
1 &   15.00 &    $-3\pi/8$  &     $7\pi/8$  &    $-7\pi/8$      & 23 &  23 &   0.0785 &   0.0812   & 15 &  15 &   0.0162  &  0.0149\\
\hline
1 &   1.05 &    $-7\pi/8$  &     $3\pi/8$  &    $-3\pi/8$      &15& 15&   0.0239  &  0.0298     &15& 15&   0.0151 &   0.0159 \\
1 &    1.10 &    $-7\pi/8$  &     $3\pi/8$  &    $-3\pi/8$      &15& 15&  0.0232  &  0.0260  &15& 15&   0.0154  &  0.0157 \\
1 &    1.50 &    $-7\pi/8$  &     $3\pi/8$  &    $-3\pi/8$      &15& 15&    0.0258  &  0.0261   &15& 15&   0.0154 &   0.0145\\
1 &    2.00 &    $-7\pi/8$  &     $3\pi/8$  &    $-3\pi/8$      &15& 15&  0.0330  &  0.0386  &15& 15&   0.0179  &  0.0158 \\
1 &   10.00 &    $-7\pi/8$  &     $3\pi/8$  &    $-3\pi/8$      &23& 23&   0.0634  &  0.0618    &15& 15&   0.0147  &  0.0157\\
1 &   15.00 &    $-7\pi/8$  &     $3\pi/8$  &    $-3\pi/8$      & 39 &  39&    0.1205  &  0.1278    & 15 &  15 &  0.0148  &  0.0150 \\
\hline
1 &    1.05 &    $-7\pi/8$  &     $7\pi/8$  &    $-7\pi/8$      &15& 15&   0.0237  &  0.0258     &15& 15&   0.0153  &  0.0159 \\
1 &    1.10 &    $-7\pi/8$  &     $7\pi/8$  &    $-7\pi/8$      &15& 15&    0.0253  &  0.0263   &15&  15&  0.0156  &  0.0164\\
1 &    1.50 &    $-7\pi/8$  &     $7\pi/8$  &    $-7\pi/8$      &15& 15&   0.0275  &  0.0334   &15&  15&  0.0160  &  0.0155\\
1 &    2.00 &    $-7\pi/8$  &     $7\pi/8$  &    $-7\pi/8$      &15& 15&    0.0311  &  0.0327    &15& 15&   0.0160  &  0.0167 \\
1 &   10.00 &    $-7\pi/8$  &     $7\pi/8$  &    $-7\pi/8$      &23& 31&  0.0597  &  0.0725     &15& 15&   0.0178  &  0.0199\\
1 &   15.00 &    $-7\pi/8$  &     $7\pi/8$  &    $-7\pi/8$      & 39 &  39 &   0.1249  &  0.1361    & 15 &   15 &   0.0148  &  0.0147 \\
\hline
\end{tabular}
\caption{Solutions that are the image of an annulus and the lower-dimensional problem with radii $a$ and $b$, and a range of inclination angles shown in sub-columns. Also shown are the number of Chebyshev points $n$ and the elapsed time in seconds.  The shaded sub-columns indicate the corresponding data.}
\label{table_ann_2}
\end{table}

\begin{table}[h!]
	\centering
	\begin{tabular}{||l|l|l|l| >{\columncolor{gray!15}}l|l|>{\columncolor{gray!15}}l|l|>{\columncolor{gray!15}}l|l|>{\columncolor{gray!15}}l|l|>{\columncolor{gray!15}}l ||}
		\hline
		\multicolumn{5}{||c|}{}& \multicolumn{4}{|c|}{Radially symmetric} & \multicolumn{4}{|c||}{Lower dimensional} \\
		\hline
		$a$ & $b$ & $\psi_a$ & \multicolumn{2}{|c|}{$\psi_b$} & \multicolumn{2}{|c|}{$n$} & \multicolumn{2}{|c|}{Elapsed time} & \multicolumn{2}{|c|}{$n$} &  \multicolumn{2}{|c||}{Elapsed time} \\
		\hline \hline
    0.1 &      1 &  $-3\pi/8$ &   $ 3\pi/8$ &   $ -3\pi/8$  & 15 & 15 &  0.0053  &  0.0038  & 15 & 15 &  0.0024  &  0.0025 \\
0.5 &    1   &  $-3\pi/8$ &    $3\pi/8$ &   $ -3\pi/8$  & 15 & 15 &  0.0048  &  0.0034  & 15 & 15 &  0.0022  &  0.0028\\
\hline
0.1 &    1  &  $-3\pi/8$ &   $ 5\pi/8$ &   $-5\pi/8$   & 15 & 15 &  0.0317  &  0.0279  & 15 & 15 &  0.0262  &  0.0282 \\
0.5 &    1  &  $-3\pi/8$ &   $ 5\pi/8$ &   $-5\pi/8$  & 15 & 15 &  0.0396  &  0.0348  & 15 & 15 &  0.0228  &  0.0228 \\
\hline
0.1 &    1  &  $-5\pi/8$ &   $ 3\pi/8$ &    $-3\pi/8$  & 15 & 15 &  0.0294 &   0.0284  & 15 & 15 &  0.0176  &  0.0168 \\
0.5 &    1  & $-5\pi/8$ &    $3\pi/8$ &    $-3\pi/8$  & 15 & 15 &  0.0247  &  0.0239  & 15 & 15 &  0.0153  &  0.0150 \\
\hline
0.1 &    1  & $-5\pi/8$ &   $ 5\pi/8$ &   $-5\pi/8$   & 15 & 15 &  0.0250  &  0.0234  & 15 & 15 &  0.0160  &  0.0163 \\
0.5 &    1 &  $-5\pi/8$ &   $ 5\pi/8$ &   $-5\pi/8$   & 15 & 15 &  0.0290  &  0.0251  & 15 & 15 &  0.0157   & 0.0156\\
\hline
\end{tabular}
\caption{Solutions that are the image of an annulus and the lower-dimensional problem with radii $a$ and $b$, and two choices of inclination angles shown in sub-columns. Also shown are the number of Chebyshev points $n$ and the elapsed time in seconds.  The shaded sub-columns indicate the corresponding data.}
\label{table_ann_3}
\end{table}

\section{{\bf P3} Study}
\label{Study3}

We conclude with our treatment of the lower-dimensional case.  We use the exact same intial guess constructions from the annular case, and this choice has robust performance in Newton's method.  The main difference between this case and what we have considered up to now is the system of differential equations is somewhat simpler, with no removable singularity, and solutions are independent of horizontal translations.  We have
\begin{verbatim}
N = @(v) [ D11*v - v(end).*cos(D03*v)
           D12*v - v(end).*sin(D03*v)
           D13*v - kappa*v(end).*D02*v
           dT1n1*v - a; dT1p1*v - b
           dT3n1*v - psia; dT3p1*v - psib ];
\end{verbatim}
and
\begin{verbatim}
L = @(v) [ D1, Z0, spdiags(v(end)*sin(D03*v),0,n-1,n-1)*D0, -cos(D03*v)
           Z0, D1, spdiags(-v(end)*cos(D03*v),0,n-1,n-1)*D0, -sin(D03*v)
           Z0, -kappa*v(end)*D0, D1, - kappa*D02*v
           dT1n1; dT1p1
           dT3n1; dT3p1 ];
\end{verbatim}

\begin{figure}[!h]
	\centering
	\scalebox{0.35}{\includegraphics{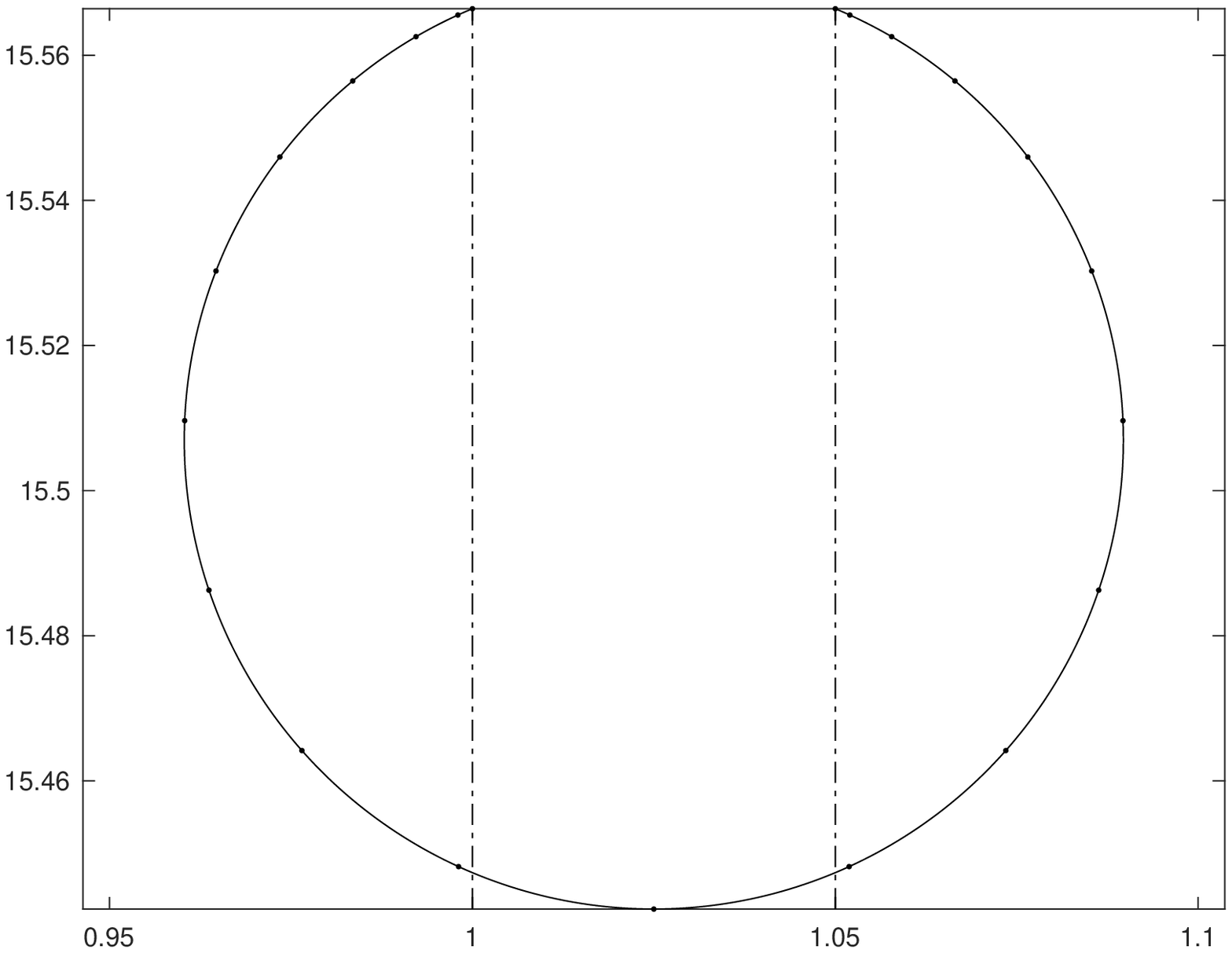}}
		\scalebox{0.35}{\includegraphics{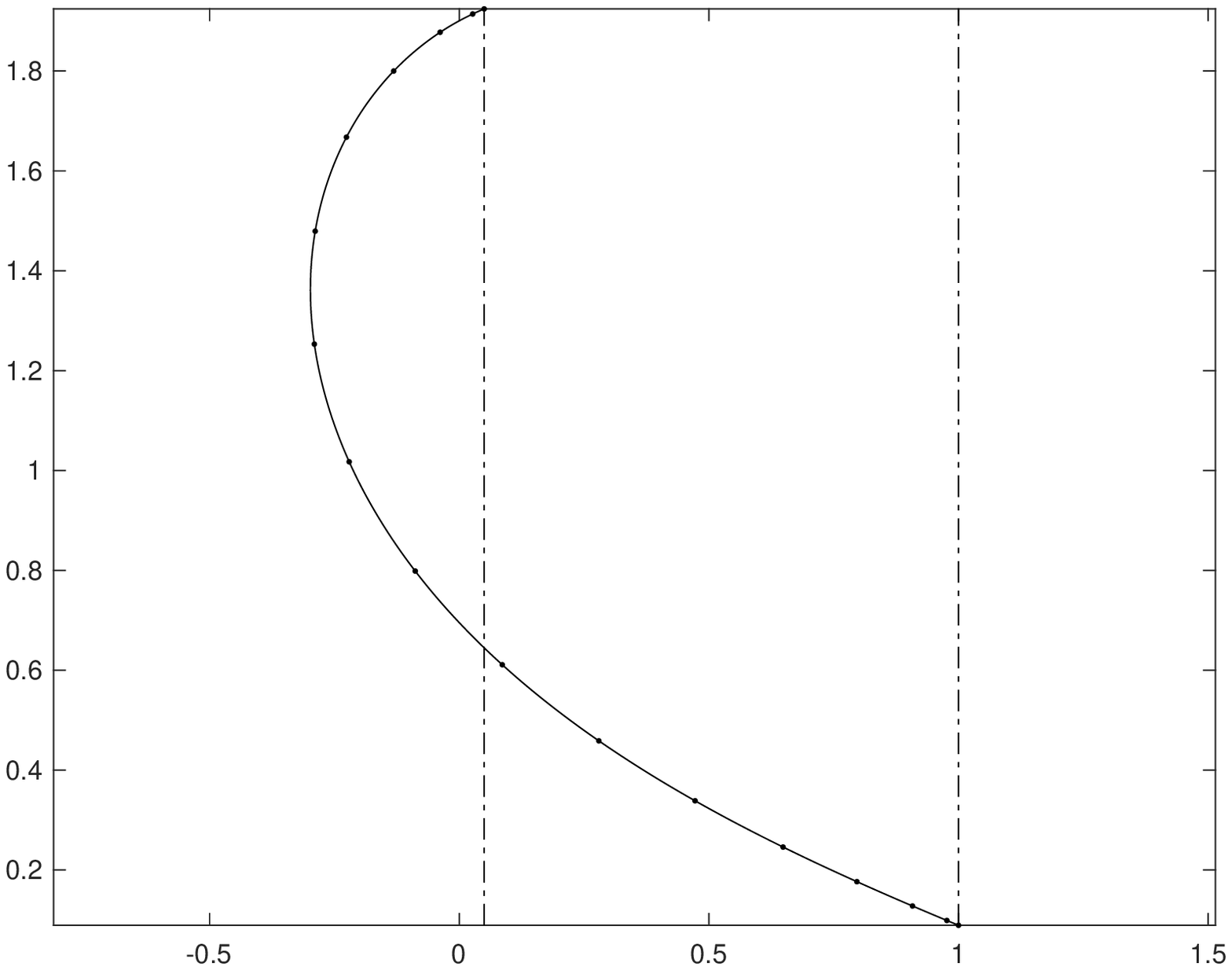}}
	\caption{A 2D surface  with inner radius $a = 1$ and outer radius $b = 1.05$ and  inclination angles $\psi_a = -7\pi/8$ and $\psi_b = 7\pi/8$ (L) and another with inner radius $a = 0.05$ and outer radius $b = 1$ and inclination  angles $\psi_a = -7\pi/8$ and $\psi_b= -\pi/8$.}
	\label{fig:2d_1_105}
\end{figure}

Then, in the absence of  that removable singularity, the solver has no numerical difficulties.  See Figure~\ref{fig:2d_1_105} (L).

We can easily compute interfaces with vertical points as in Figure~\ref{fig:2d_1_105} (R).  To contrast the behavior of this system with the rotationally symmetric systems, the example given here has $0 < a < b$, but the solution curve itself crosses the vertical axis.  This cannot happen in the previous problem, which is a consequence of the regularity theory for the full dimensional and radially symmetric problems.

\begin{figure}[!h]
	\centering
	\scalebox{0.35}{\includegraphics{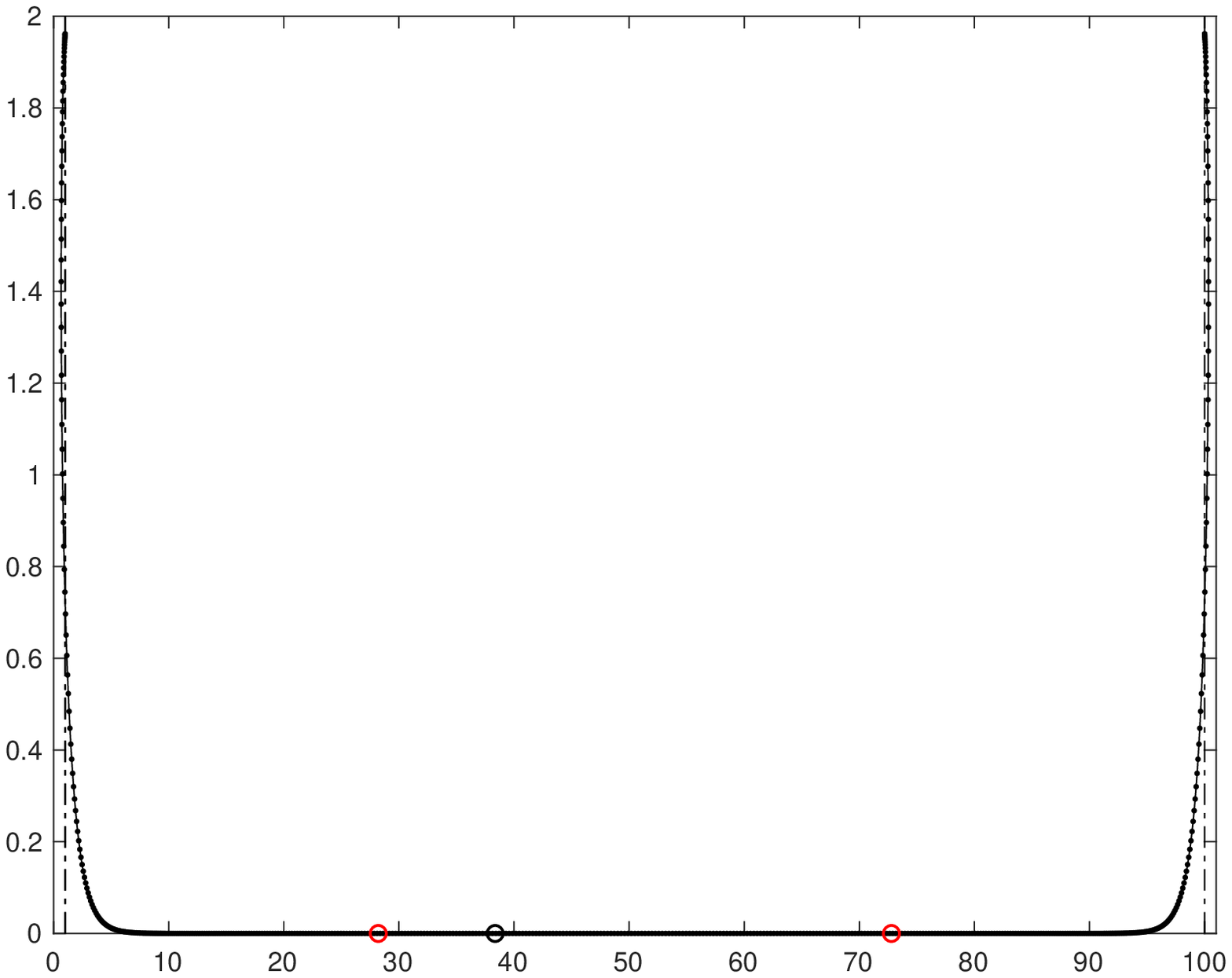}}
	\scalebox{0.35}{\includegraphics{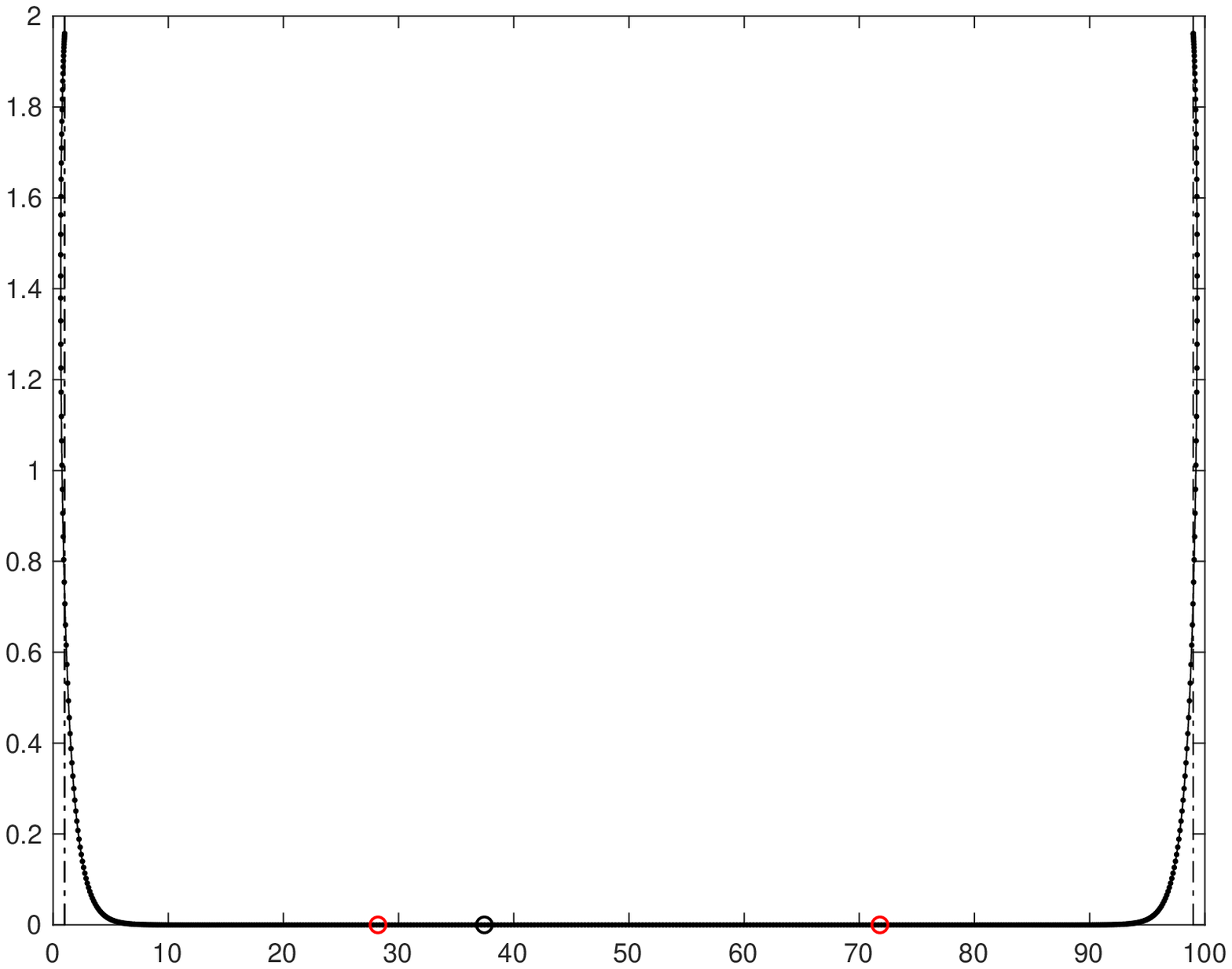}}
	\caption{On the left is a 2D surface  with inner radius $a = 1$ and outer radius $b = 100$ and  inclination angles $\psi_a = -7\pi/8$ and $\psi_b = 7\pi/8$ using 541 data points along the interface.  On the right $b = 99$, and the other parameters remain the same, with 541 data points also being used.  The axes are scaled differently for visibility, and the minimum point on each curve is circled.}
	\label{fig:2d_big}
\end{figure}

Finally, we explore the solver's performance as the gap between $a$ and $b$ gets somewhat large.  On the left of Figure~\ref{fig:2d_big} we take $a = 1$ and $b = 100$ and we are able to obtain the requested precision with 541 Chebyshev points.  We did change the default values of the loop counters to  $\texttt{max\_iter\_newton = 100}$ and $\texttt{max\_iter\_bvp = 500}$.  If we had used the shooting method sketched in the previous section, the nonlinear solver would not have been able to produce a horizontal point $(m,u_m)$ to meet the boundary conditions within any reasonable accuracy, as the most minute variations of $u_m$ lead to wild variations of the interface when $u_m$ is this close to zero.  On the right of   Figure~\ref{fig:2d_big}, $b$ was changed to 99, and the same number of data points were used in the adaptive algorithm.  To compare the two figures, for  $b = 100$ the minimum height $u_m  = $ \texttt{-1.6098e-15} at $\tau  = $ \texttt{-0.2375} and $m =  $ \texttt{38.3626} there.  For $b = 99$, $u_m = $ \texttt{-1.5543e-15} at $\tau = $\texttt{-0.2481} and we have $m = $ \texttt{37.4709}  
with.  In both cases, the minimum value is within the tolerance of 0, and about an order of magnitude away from machine epsilon.  The intersections of both curves with the height $\texttt{tol\_bvp}$ was determined.  The intersection points are circled in red. For $b = 100$, the first intersection was at $R =       28.2080$ and the second at $R =  72.7924$, for a total length of 44.5844 along the $R$ axis where the solution curve had a height within tolerance of zero.  Here a notational convention has been used that streamlines the different codes written.  We use $R = X$ and $r = x$ in the code to avoid internal variable conflicts.  This abuse of notation is used here as well.  Then, for $b = 99$, the first intersection was at $R =    28.2099$ and the second at $R = 71.7915$ for a total length of 43.5816 along the $R$ axis where the solution curve had a height within tolerance of zero.

\section{Conclusions and closing remarks}
\label{conclusions}

We have presented Chebyshev spectral methods for three capillary surfaces with a variety of prescribed inclination angles that are met at prescribed radii.  For modest problems, the convergence to solutions is extremely fast and uses very little memory.

For more challenging problems the adaptive algorithm automatically increases the number of Chebyshev points to achieve the prescribed tolerances.  This precess has worked very well in almost all of the cases we have found.  It is possible sometimes to break these codes for extremes of the inclination angles near $\pm\pi$ and either radii too close to each other with $\psi_a\psi_b > 0$, inner radius $a$ too close to 0, $b$ too large, or the  radii too far apart .  Typically this happens when $|\psi_a|, |\psi_b| > \pi/2$, and the closer to $\pm\pi$, the more challenging.   For these multi-scale problems customizing the initial guess and carefully tuning the increase of the Chebyshev nodes inside the adaptive part of the algorithm can lead to success when the base code does not converge.  Still, some rather extreme problems may  not work even then, and for those problems a multi-scale approach is needed.  This is the subject of a future paper.

These codes are prototypes, and can be adapted to similar problems in fluid mechanics.  For our exposition, we have included $|\psi_b | > \pi/2$, which does not immediately appear as a physical surface in any of the physical problems we discussed in the introduction.  There are other geometric constraints that one might put on the equilibrium shapes of fluids to model different fluid problems, and we have only applied these methods to a few of the possible examples.


\end{document}